\documentclass[10pt,leqno,english]{article}

 \usepackage{latexsym}
 \usepackage[T1]{fontenc}
 \usepackage{amssymb}
 \usepackage{amsmath}
 \usepackage{amsthm}
\usepackage{natbib}

 \usepackage[dvips]{graphics}
\usepackage{graphicx}

 \newcommand{\Er}{\mathcal{O}}

 \newcommand{ \un }{\mathbb{I} }
 \newcommand{ \p }{\mathbb{P} }
 \newcommand{ \pa }{\mathbb{P}^{\alpha} }
  
\newcommand{ \pam }{\mathbb{P}^{\alpha}_{\tmo}}
 \newcommand{ \E }{\mathbb{E}}
 \newcommand{ \Ea }{\mathbb{E}^{\alpha}}
 \newcommand{ \Eam }{\mathbb{E}^{\alpha}_{\tmo}}

 \newcommand{ \R }{ \mathbb{R} }
 \newcommand{ \Z }{ \mathbb{Z} }
 \newcommand{ \N }{ \mathbb{N} }

 \newcommand{ \rw}{\textrm{R.W.R.E. }}

 \newcommand{ \f }{ \mathcal{F} }

 \newcommand{ \tm }{ \tilde{m} }
 \newcommand{ \tM }{ \tilde{M} }
\newcommand{ \tmo }{ \tilde{m}_0 }

\newcommand{ \G}{ \mathcal{G} }

\newcommand{ \A}{ \mathcal{A} }
 
 \newcommand{ \M }{ \mathcal{M} }

\newtheorem{The}{{\bf Theorem}}[section]

\theoremstyle{definition}
\newtheorem{Def}[The]{{\bf Definition}}
\theoremstyle{plain}
 \newtheorem{Lem}[The]{Lemma}
 \newtheorem{Cor}[The]{\bf Corollary}
 \newtheorem{Pro}[The]{\bf Proposition}
 \theoremstyle{definition}
\newtheorem{Rem}[The]{{\bf Remark}}

 \newenvironment{Pre}{\noindent \textbf{Proof.} \\ }{$\
 \blacksquare$}

\newenvironment{Prepr}[1]{\noindent \textbf{Proof (of Proposition #1).} \\ }{$\
 \blacksquare$}

\makeatletter

\@addtoreset{equation}{section} \makeatother

\title{ Alternative proof for the localization of Sinai's walk
\\ \vspace{1cm}
 \large{Pierre Andreoletti} $^{\dag}$,\footnote{ Centre de Physique Théorique, C.N.R.S. UMR 6207, Université Aix-Marseille I, II, Université du Sud-Toulon-Var, F.R.U.M.A.M. (Marseille, France)
 and Centro de Modelamiento Mathematico C.N.R.S. U.M.R 2071, Universidad de Chile (Santiago,
 Chile). \newline \vspace{0.1cm}  $\quad$  MSC 2000 60G50.  }}

\begin{document}

\bibliographystyle{unsrtnat}

\maketitle

\noindent  $^{\dag}$ Université Aix-Marseille II, Facult\'e des
sciences de Luminy, C.P.T. case 907, 13288 Marseille cedex 09
France. e-mail  \texttt{andreole@cpt.univ-mrs.fr}

\noindent \\ \textbf{Abstract:} We give an alternative proof of
the localization of Sinai's random walk in random environment
under weaker hypothesis than the ones used by Sinai. Moreover we
give estimates that are stronger than the one of Sinai on the
localization neighborhood and on the probability for the random
walk to stay inside this neighborhood.

 \noindent \\  \textit{Key
words and phrases : Random environment, random walk, Sinai's
regime, Markov chain.}

\noindent \\
CPT-2004/P.068

\section{Introduction}

Random Walks in Random Environment (R.W.R.E.) are basic processes
in random media. The one dimensional case with nearest neighbor
jumps, introduced by \cite{Solomon}, was first studied by
\cite{KesKozSpi}, \cite{Sinai}, \cite{Golosov},  \cite{Golosov0}
and \cite{Kesten2} all these works show the diversity of the
possible behaviors of such walks depending on hypothesis assumed
for the environment. At the end of the eighties \cite{Deh&Revesz}
and \cite{Revesz} give the first almost sure behavior of the \rw
in the recurrent case. Then we have to wait until the middle of
the nineties to see new results. An important part of these new
results concerns the problem of large deviations first studied by
\cite{GreHol2} and then by \cite{ZeiGan}, \cite{PiPO},
\cite{PiPoZe} and \cite{ZeCoGa} (see \cite{Zeitouni} for a
review). In the same period using the stochastic calculus for the
recurrent case \cite{Shi},
 \cite{HuShi2}, \cite{HuShi1}, \cite{Hu1}, \cite{Hu} and \cite{HuShi0} follow the works of
\cite{Schumacher} and \cite{Brox}  to give very precise results on
the random walk and its local time (see \cite{Shi1} for an
introduction). Moreover recent results on the problem of aging are
given in \cite{ZeDEGu}, on the moderate deviations in
\cite{ComPop} for the recurrent case, and on the local time in
\cite{GanShi} for the transient case. In parallel to all these
results a continuous time model has been studied, see for example
\cite{Schumacher} and \cite{Brox}, the works of \cite{Tanaka},
\cite{Mathieu2}, \cite{Tanaka2}, \cite{KawTan}, \cite{Mathieu1}
and \cite{Taleb0}.

Since the beginning of the eighties the delicate case of \rw in
dimension larger than 2 has been studied a lot, see for example
\cite{Kalikow}, \cite{AnKhSi}, \cite{Durett}, \cite{BoCoGeDo}, and
\cite{BricKup}. For recent reviews (before 2002) on this topics
see the papers of \cite{ Sznitman} and \cite{Zeitouni}. See also
\cite{Sznitman3}, \cite{Varadhan}, \cite{Rass2} and
\cite{ZeiCom2}.

In this paper we are interested in Sinai's walk i.e the one
dimensional random walk in random environment with three
conditions on the random environment: two necessaries hypothesis
to get a recurrent process (see \cite{Solomon}) which is not a
simple random walk and an hypothesis of regularity which allows us
to have a good control on the fluctuations of the random
environment.

The asymptotic behavior of such walk was discovered by
\cite{Sinai}, he showed that this process is sub-diffusive and
that at time $n$ it is localized in the neighborhood of a well
defined point of the lattice. This \textit{point of localization}
 is a random variable depending only on the random environment
and $n$, its explicit limit distribution was given, independently,
by \cite{Kesten2} and \cite{Golosov0}.

\noindent Here we give an alternative proof of Sinai's results
under a weaker hypothesis. First we recall an elementary method
proving that for a given instant $n$ Sinai's walk is trapped in a
basic valley denoted $\{\tM_0',\tmo,\tM_0\}$ depending only on $n$
and on a realization of the environment. Then we give a proof of
the localization, this proof is based on an analysis of the return
time to $\tmo$. We get a stronger result than Sinai : we find that
a size of the neighborhood of the localization depends on $n$ like
$(\log_2 n)^{9/2} (\log n)^{3/2}$ instead of $\delta (\log n)^2$
found by Sinai. Moreover we compute the rates of the convergence
of the probabilities (for the random walk and the random
environment). Our method is based on the classification of the
valleys obtained by ordered refinement of the basic valley
$\{\tM_0',\tmo,\tM_0\}$. The properties of the valleys obtained by
this operation are proved with some details.

\noindent \\ This paper is organized as follows. In section 2 we
describe the model, we give some basic notions on the random
environment and present the main results.
 In section \ref{sec4}
 we give the properties of the random
environment needed in section \ref{sec5} to prove the main
 results. In the Appendix we make the proof of the properties of the random environment.


\section{Description of the model and main results}

\subsection{Sinai's random walk definition }

Let $\alpha \equiv (\alpha_i,i\in \Z)$ be a sequence of i.i.d.
random variables taking values in $(0,1)$ defined on the
probability space $(\Omega_1,\f_1,Q)$, this sequence will be
called random environment. A random walk in random environment
(denoted R.W.R.E.) $(X_n,n
\in\N)$ is a sequence of random variable taking value in $\Z$, defined on $( \Omega,\f,\p)$ such that \\
$ \bullet $ for every fixed environment $\alpha$, $(X_n,n\in \N)$
 is a Markov chain with the following transition probabilities, for all
 $i\in \Z$
 \begin{eqnarray}
& & \p^{\alpha}\left[X_n=i+1|X_{n-1}=i\right]=\alpha_i, \label{mt} \\
& & \p^{\alpha}\left[X_n=i-1|X_{n-1}=i\right]=1-\alpha_i \equiv
\beta_i. \nonumber
\end{eqnarray}
We denote by $(\Omega_2,\f_2,\pa)$ the probability space
associated to this Markov chain. \\
 $\bullet$ $\Omega = \Omega_1 \times \Omega_2$, $\forall A_1 \in \f_1$ and $\forall A_2 \in \f_2$,
$\p\left[A_1\times
A_2\right]=\int_{A_1}Q(dw_1)\int_{A_2}\p^{\alpha(w_1)}(dw_2)$.

\noindent \\ The probability measure $\pa\left[\left.
.\right|X_0=a \right]$  will be  denoted $\pa_a\left[.\right]$,
 the expectation associated to $\pa_a$: $\Ea_a$, and the expectation associated to
 $Q$:
 $\E_Q$.

\noindent \\ Now we introduce the hypothesis we use in all this
work. Denoting $\left(\epsilon_i=\log [(1-\alpha_i)/\alpha_i],
i\in \Z \right)$, the two following hypothesis are the necessaries
hypothesis
\begin{eqnarray}
 \E_Q\left[\epsilon_0 \right]=0 , \label{hyp1bb} \label{hyp1}
\end{eqnarray}
\begin{eqnarray}
\E_Q\left[\epsilon_0^2\right]\equiv \sigma^2
>0 . \label{hyp0}
\end{eqnarray}
 \cite{Solomon} shows that under \ref{hyp1} the process $(X_n,n\in \N)$
is $\p$ almost surely null recurrent and \ref{hyp0} implies that
the model is not reduced to the simple random walk. In addition to
\ref{hyp1} and \ref{hyp0} we will consider the following
hypothesis of regularity, there exists $\kappa^+ \in \R^{*}_+$
such that for all $\kappa \in ]0,\kappa^+[$
\begin{eqnarray}
\E_Q\left[e^{\kappa \epsilon_0}\right]<\infty \textrm{ and } \
\E_Q\left[e^{-\kappa \epsilon_0}\right]<\infty.  \label{hyp4}
\end{eqnarray}

\noindent We call \textit{Sinai's random walk} the random walk in
random environment previously defined with the three hypothesis
\ref{hyp1}, \ref{hyp0} and \ref{hyp4}.

\noindent  Notice that Y. Sinai used the stronger hypothesis :
\begin{eqnarray}
\alpha_0\geq \textrm{const} >0,\  1-\alpha_0 \geq \textrm{const}
>0. \label{hyp2}
\end{eqnarray}



\subsubsection{The random potential and the valleys }

\begin{Def} \label{defpot2} The random potential $(S_k,\  k \in
\R)$ associated to the random environment $\alpha$ is defined by
\begin{eqnarray}
 S_k=\left\{ \begin{array}{ll} \sum_{1\leq i \leq k} \epsilon_i, &  k=1,2,\cdots , \\
   \sum_{k \leq i \leq -1} \epsilon_i , &  k=-1,-2,\cdots  , \end{array} \right
. \nonumber
 \end{eqnarray}
 for the other $k\in \R\!\smallsetminus\!\Z $, $(S_k,k)$ is
defined by linear interpolation,  and $S_0=0$. We denote $(S^n_t,t
\in \R)$ the normalized potential
 associated to $(S_k,\  k \in
\Z)$
\begin{eqnarray}
 S_k^n= \frac{S_k}{ \log n },\ k\in \Z .
 \end{eqnarray}
\end{Def}

\begin{Def} \label{c2s2d1}
 We will say that the triplet $\{\tilde{M}',\tilde{m},\tilde{M}''\}$ is a \textit{valley} if
 \begin{eqnarray}
 S^{n}_{\tilde{M}'}=\max_{\tilde{M}' \leq t \leq \tilde{m}} S^{n}_t ,  \\
 S^{n}_{\tilde{M}''}=\max_{\tilde{m} \leq t \leq
 \tilde{M''}}S^{n}_t ,\\
S^{n}_{\tilde{m}}=\min_{\tilde{M}' \leq t \leq \tilde{M}''}S^{n}_t
\ \label{2eq58}.
 \end{eqnarray}
If $\tilde{m}$ is not unique, we choose the one with the smallest
absolute value.
 \end{Def}

\begin{Def} \label{deprofvalb}
 We will call \textit{depth of the valley} $\{\tM',\tm,\tM''\}$ and we
 will denote $d([\tilde{M}',\tilde{M}''])$ the quantity
\begin{eqnarray}
 \min(S^{n}_{\tilde{M}'}-S^{n}_{\tilde{m}},S^{n}_{\tilde{M}''}-S^{n}_{\tilde{m}})
 .
 \end{eqnarray}
 \end{Def}

\noindent  Now we define the operation of refinement.
 \begin{Def}
Let  $\{\tM',\tm,\tM''\}$ be a valley. Let
  $\tM_1$ and $\tm_1$ be such that $\tm \leq \tM_1< \tm_1 \leq \tM''$
  and
 \begin{eqnarray}
 S_{\tM_1}^n-S_{\tm_1}^n=\max_{\tm \leq t' \leq t'' \leq
 \tM''}(S_{t'}^n-S_{t''}^n) .
 \end{eqnarray}
 We say that the couple $(\tM_1,\tm_1)$ is obtained by a \textit{right refinement} of $\{\tM',\tm,\tM''\}$. If the couple $(\tm_1,\tM_1)$ is not
 unique, we will take  the ones such that $\tm_1$ and $\tM_1$ have the smallest  absolute value. In a similar way we
  define the left refinement operation.
 \end{Def}

\noindent In all this work, we denote $\log_p $ with $p\geq 2$ the
$p$ iterated logarithm and we assume that $n$ is large enough such
that $\log_p n$ is positive. Let $\gamma>0$ a free parameter,
denoting $\gamma(n)= (\gamma \log_2 n) (\log n)^{-1} $ we  define
what we will call a \textit{valley containing $0$ and of depth
larger than $1+\gamma(n)$}.

\begin{Def} \label{thdefval1b} \label{thdefval1} For $\gamma>0$ and $n>3$, we say that a valley $\{\tM',\tm,\tM''\}$ contains $0$ and
is of depth larger than $1+ \gamma(n) $ if and only if
\begin{enumerate}
\item  $ 0 \in [\tM',\tM'']$, \item $d\left(\{\tM',\tM''\}\right)
\geq 1+\gamma(n) $ , \item if $\tm<0,\ S_{\tM''}^n-\max_{ \tm \leq
t \leq 0
}\left(S_t^n\right) \geq \gamma(n) $ , \\
 if $\tm>0,\
S_{\tM'}^n-\max_{0 \leq t \leq
 \tm}\left(S_t^n\right) \geq \gamma(n) $ .
\end{enumerate}
\end{Def}

\subsubsection{The basic valley $\{\tM_0',\tmo,\tM_0\}$ \label{2.1.2}}
\noindent  We recall the notion of \textit{basic valley},
introduced by Y. Sinai and denoted here $\{\tM_0',\tmo,\tM_0\}$.
The definition we give is inspired by the work of \cite{Kesten2}.
First let $\{\tM',\tmo,\tM''\}$ be the smallest valley that
contains $0$ and of depth larger than $1+\gamma(n)$. Here smallest
means that if we construct, with the operation of refinement,
other valleys in $\{\tM',\tmo,\tM''\}$ such valleys will not
satisfy one of  the properties of Definition \ref{thdefval1b}.
$\tM_0'$ and
$\tM_0$ are defined from $\tmo$ in the following way \\
if $\tmo>0$
\begin{eqnarray}
& & \tM_0'=\sup \left\{l\in \Z_-,\ l<\tmo,\ S_l^n-S_{\tmo}^n\geq
1+\gamma(n),\
S_{l}^n-\max_{0 \leq k \leq \tmo}S_k^n \geq \gamma(n) \right\} ,\\
& & \tM_0=\inf \left\{l\in \Z_+,\ l>\tmo,\ S_l^n-S_{\tmo}^n\geq
1+\gamma(n)\right\} . \label{4.8}
\end{eqnarray}
If $\tmo<0$
\begin{eqnarray}
& & \tM_0'=\sup \left\{l\in \Z_-,\ l<\tmo,\ S_l^n-S_{\tmo}^n\geq
1+\gamma(n)\right\} , \\
& & \tM_0=\inf \left\{l\in \Z_+,\ l>\tmo,\ S_l^n-S_{\tmo}^n\geq
1+\gamma(n),\ S_{l}^n-\max_{ \tmo \leq k \leq 0}S_k^n \geq \gamma(
n) \right\} . \label{4.10}
\end{eqnarray}
If $\tmo=0$
\begin{eqnarray}
& & \tM_0'=\sup \left\{l\in \Z_-,\ l<0,\ S_l^n-S_{\tmo}^n\geq
1+\gamma(n) \right\} , \\
& &  \tM_0=\inf \left\{l\in \Z_+,\ l>0,\ S_l^n-S_{\tmo}^n\geq
1+\gamma(n) \right\} . \label{4.12}
\end{eqnarray}

\noindent One can ask himself if the basic valley exists, in the
Appendix \ref{sec6} we prove the following lemma :

\begin{Lem} \label{moexiste} Assume \ref{hyp1}, \ref{hyp0} and \ref{hyp4}, for all $\gamma>0$ there exists $n_0\equiv n_0(\gamma,\sigma,E[|\epsilon_0|^3])$
such that for all $n>n_0$
\begin{eqnarray}
Q\left[\{\tM_0',\tmo,\tM_0\} \neq \varnothing \right] \geq 1- (6
\gamma \log_2 n)(\log n)^{-1}.
\end{eqnarray}
\end{Lem}

\begin{Rem}
\noindent In all this paper we use the same notation $n_0$ for an
integer that could change from line to line. Moreover in the rest
of the paper we do not always make explicit the dependance on
$\gamma$ of all those $n_0$ even if Lemma \ref{moexiste} is
constantly used.
\end{Rem}

\subsection{Main results : localization phenomena}

The following result shows that Sinai's random walk is
sub-diffusive :

\begin{Pro}\label{lem2} There exists a strictly positive numerical constant $h>0$, such that if  \ref{hyp1} and \ref{hyp0} hold and for all $\kappa \in ]0,\kappa^+[$ \ref{hyp4} hold, for all  $\gamma>2$ there exists $n_0\equiv
n_0(\gamma)$ such that for all $n>n_0$, there exists $G_n \subset
\Omega_1$ with $Q\left[G_n\right] \geq 1- h \left((\log_3 n)(
\log_2 n)^{-1} \right)^{1/2}$ and
\begin{eqnarray}
\sup_{\alpha \in G_n} \left\{
\p^{\alpha}_0\left[\bigcup_{m=0}^{n}\left\{X_m\notin\left[\tM_0',\tM_0
\right]\right\}\right] \right\} \leq
 \frac{2 \log_2 n}{\sigma^2(\log n)^{\gamma-2}} , \label{2.19}
 \end{eqnarray}
 moreover
\begin{eqnarray}
\sup_{\alpha \in G_n} \left\{
\p^{\alpha}_0\left[\bigcup_{m=0}^{n}\left\{X_m\notin\left[-(\sigma^{-1}
\log n)^2 \log_2 n, (\sigma^{-1} \log n)^2 \log_2 n
\right]\right\}\right] \right\} \leq
 \frac{2 \log_2 n}{\sigma^2(\log n)^{\gamma-2}} . \label{2.20}
 \end{eqnarray}
\end{Pro}

\begin{Rem}
A weaker form of this result can be found in the paper of
\cite{Sinai} (Lemma 3 page 261). The set $G_n$ is called set of
"good" environments. We will define it precisely in section
\ref{sec4}. This set is defined by collecting all the properties
on the environment we need to prove
our results. \\
 \ref{2.19} shows that Sinai's walk is
trapped in the basic valley $\{\tM_0',\tmo,\tM_0\}$ which is
random,  depending only on the random media and on $n$. More
precisely, using \ref{2.20}, with an overwhelming probability
$\{\tM_0',\tmo,\tM_0\}$ is within an interval centered at the
origin and of size $2(\sigma^{-1} \log n)^2 \log_2 n $. \noindent
In all this work $h$ is a strictly positive numerical constant
that can grow from line to line if needed.
\end{Rem}

\noindent The following remarkable result was proved by
 \cite{Sinai}

\begin{The}
Assume \ref{hyp1}, \ref{hyp0} and \ref{hyp2}, for all $\epsilon>0$
and all $\delta>0$ there exists $n_0\equiv n_0(\epsilon,\delta)$
such that for all $n>n_0$, there exists $C_n \subset \Omega_1$
with $Q\left[C_n\right] \geq 1- \epsilon $ and
\begin{eqnarray}
\lim_{n\rightarrow +\infty } \sup_{\alpha \in G_n}
\pa_0\left[\left|\frac{X_n}{\log^2n}-m_{0}\right|
> \delta \right]=0,
\end{eqnarray}
$m_0= \tmo (\log n)^{-2}$.
\end{The}

\noindent In this paper we improve Sinai's result in the following
way, for all $\kappa \in ]0,\kappa^+[$ we denote
$\gamma_0=\frac{12}{\kappa}+\frac{21}{2}$,

\begin{The} \label{thSinai1}
 There exists a strictly positive numerical constant $h>0$, such that if  \ref{hyp1} and \ref{hyp0} hold and for all $\kappa \in ]0,\kappa^+[$ \ref{hyp4} hold, for all $\gamma>\gamma_0$ there exists $n_0\equiv
n_0(\gamma)$ such that for all $n>n_0$, there exists $G_n \subset
\Omega_1$ with $Q\left[G_n\right] \geq 1- h \left((\log_3 n) (
\log_2 n)^{-1} \right)^{1/2}$ and
\begin{eqnarray}
\sup_{\alpha \in G_n}
\left\{\p^{\alpha}_0\left[\left|\frac{X_n}{\log^2n}-m_{0}\right|
> \G \gamma \frac{  (\log_2 n )^{9/2}}{(\log n)^{1/2}}
\right] \right\}\leq \frac{4(\log_2 n)^{9/2}}{\sigma^{10}(\gamma
\log n )^{\gamma-\gamma_0}} , \label{eqthSinai21b}
\end{eqnarray}
$m_0= \tmo (\log n)^{-2}$ and $\G= (1600)^2$.
\end{The}

\begin{Rem}
\noindent This result shows that, for a given instant $n$
sufficiently large, with a $Q$ probability tending to one,
 $X_n$ belongs to a neighborhood of the point
$\tmo $ with a $\pa$ probability tending to one. The size of this
neighborhood is of order $(\log n)^{3/2} (\log_2 n)^{9/2}$ that is
negligible comparing to the typical range of Sinai's walk of order
$(\log n)^2$. Moreover an estimate on the rates of the convergence
of these probabilities are given but we did not try any attempts
to optimize these rates. However if we look for an annealed
result, that means a result in $\p$ probability, we get
\begin{eqnarray}
\p\left[\left|\frac{X_n}{\log^2n}-m_{0}\right| > \G \gamma \frac{
(\log_2 n )^{9/2}}{(\log n)^{1/2}} \right]\leq 2 h \left(\frac{
\log_3  n}{ \log_2 n}\right)^{1/2} \label{eqthSinai21b}
\end{eqnarray}
and the rate in $ (\log_3  n) (\log_2 n)^{-1}$ cannot be improved
to something like $(\log n)^{-a}$ with $a>0$ without changing the
size of the localization neighborhood.
 \noindent \\ We recall that the
explicit limit distribution of $m_0$ was given independently by
\cite{Kesten2} and \cite{Golosov0}. \\
\end{Rem}

\subsection{Ideas of the proofs}

In this section we describe in detail the structure of the paper
and give the main ideas of the proofs of Propositions \ref{lem2}
and Theorem \ref{thSinai1}. For these proofs we need both
\textit{arguments on the random environment} and \textit{arguments
on
the random walk}. \\
Because of the technical aspect of the arguments on the
environment, we summarize the needed \textit{results on the
environment} in section \ref{sec4} and we have put the proofs of
these results in the Appendix at the end of the paper. So assuming
 the results of section \ref{sec4}, the proofs of the main results are
 limited to the arguments for the walk
 given in section \ref{sec5}.

\noindent \\ \textit{Results on the random environment (section
\ref{sec4})} First we describe the ordered chopping in valleys.
According to this construction, based on the refinement operation,
we get a set of valleys with the two following main properties :
1. the valleys of this set are ordered (in the sense of the depth)
2. the depth of these valleys decrease when they get close to
$\tmo$. This construction is one of the important point to get
estimations more precise than Sinai's ones, for the environment,
and therefore for the walk. We have collected all the needed
properties of the valleys in a definition (Definition
\ref{super}). All the environments that satisfy this definition
are called \textit{good environment} and we get the set of good
environment (called $G_n$, $n$ is the time). The longest part of
this work will be to prove that $Q[G_n]$ satisfies the mentioned
estimate, this is the purpose of the Appendix.

 \noindent \\ \textit{Arguments for the walk (section
\ref{sec5})}

First we recall basic results on birth and death processes used
all over the different proofs. We will always assume that the
random environments belong to the set of good environments.

 The proof of Proposition \ref{lem2} is based on
a basic argument: with an overwhelming probability, first the walk
reach the bottom of the basic valley $\tmo$ and then prefer
returning $n$ times to this point instead of climbing until the
top of the valley (\textit{i.e} reaching one of the points
$\tM_0'$ or $\tM_0$). Moreover, according to one of the properties
of the good environments, the size of the basic valley
$max\{|\tM_0'|,|\tM_0|\} \leq (\sigma^{-1} \log n)^2 \log_2 n$. So
we get the Proposition. We will see that to get this result we
have used very few properties of the good environments.

 The proof of Theorem \ref{thSinai1} is based on the two following
 facts : \textbf{Fact 1} With an overwhelming probability, the last
 return to $\tmo$ before the instant $n$, occurs at an instant
 larger than $n-q_n$. $q_n$ is a function of $n$ given by $\log q_n \thickapprox ((\log
 n)^{3/2} (\log_2 n)^{7/2})^{1/2}$.
 \textbf{Fact 2} We use the same argument of the proof of
 Proposition \ref{lem2}. With an overwhelming probability,
 starting from $\tmo$ with an amount of time $n-(n-q_n)=q_n$ the
 walk is trapped in a valley of size of order $(\log
 q_n)^2 \log_2 q_n \thickapprox (\log
 n)^{3/2} (\log_2 n)^{9/2} $. This gives the Theorem.\\
The hardest part is to prove Fact 1, for this we use both an
analysis of the return time to $\tmo$ (section \ref{3par2}) and
the ordered chopping in valleys. The main idea is to prove that
for each scale of time larger than $q_n$, the walk will return to
$\tmo$ with an overwhelming probability. These scales of time are
chosen as function of the depth of the ordered
 valleys, \textit{i.e} for each scale of time corresponds a
 valleys. What we prove is that for each scale of time the walk
 can't be trapped in the corresponding valley. Indeed, starting from $\tmo$,
 if the walk has enough time to reach the bottom of a valley  it has enough time to escape from
 it and therefore to return to $\tmo$.

\noindent \\
 \textit{Arguments for the random environment (Appendix)} While
 the proof of the results for the random environment are
technical we give some details. This provide completeness to the
present paper and shows the difficulties to work with the
hypothesis \ref{hyp4}.

\section{Good properties of a random environment \label{sec4}}

In this section we present different notions for the environment
that are used to prove the main results. We give a method to
classify some valleys obtained from $\{\tM_0',\tmo,\tM_0\}$ by the
operation of refinement. To do this we need some basic result on
$\{\tM_0',\tmo,\tM_0\}$. Then we define the set of the "good"
environments, this set contains all the environments that satisfy
the needed properties to prove the main results.

\subsection{Ordered chopping in valleys \label{OC} }

\begin{Pro} \label{8eq18} There exists $h>0$ such that if \ref{hyp1bb}, \ref{hyp0} and \ref{hyp4} hold, for all $\gamma>0$ there exists $n_0 \equiv
n_0\left(\gamma\right)$ such that for all $n>n_0$, we have
\begin{eqnarray}
& & Q\left[\tM_0 \leq  (\sigma^{-1} \log n)^2 \log_2 n \right]
\geq 1-h \left((\log_3  n) (\log_2 n)^{-1}\right)^{1/2}
 , \label{2eq126} \\
& &Q\left[\tM_0' \geq - (\sigma^{-1} \log n)^2 \log_2 n \right]
\geq 1-h \left((\log_3  n) (\log_2 n)^{-1}\right)^{1/2}.
  \label{2eq126b}
\end{eqnarray}
\end{Pro}

\noindent Before making a classification of the valleys we need to
introduce the following notations, let $\gamma>0$ and $n>3$
\begin{eqnarray}
& & b_n=[ (\gamma)^{1/2} (\log n \log_2 n)^{3/2} ], \label{4.10b}  \\
& & k_n= ((\sigma^{-1} \log n)^2 \log_2 n) (b_n)^{-1},
\label{4.11b}
\end{eqnarray}
where $[a]$ is the integer part of $a \in \R$. Using \ref{4.10b}
 and \ref{4.11b} we construct a deterministic chopping of the
interval $(-(\sigma^{-1}\log n)^2 \log_2 n,(\sigma^{-1}\log n)^2
\log_2 n) $ into pieces of length $b_n$. Moreover we define :
\begin{eqnarray}
l_n= D \sigma^2 \log k_n,\ D=1000. \label{4.12b}
\end{eqnarray}

\noindent
 We make the following construction, let us take
 $\{\tM_0',\tm_0,\tM_0\}$ as the initial valley (see Section \ref{2.1.2}). Let us denote
 $\M_0'=\{\tM_0',\tm_0\}$ and  $\M_0=\{\tm_0,\tM_0\}$. \\
First we consider the first \textit{right refinement} of the
valley
 $\{\tM_0',\tm_0,\tM_0\}$ we denote $\{\tM_1,\tm_1\}$ the couple of maximizer and
minimizer obtained after this refinement, let us add this points
to the set $\M_0$ to get $\M_0=\{\tm_0,\tM_1,\tm_1,\tM_0\}$. Now
we consider the first refinement of $\{\tm_0,\tM_1\}$, we get the
couple $\{\tM_2,\tm_2\}$ that we add to the set $\M_0$ and so on
until we obtain the points $\{\tM_r,\tm_r\}$ such that
$\tM_{r-1}-\tm_{0} \geq l_n b_n $ and $\tM_r-\tm_0 \leq l_n b_n$.
From this construction (see Figure \ref{thfignew1}) we obtain a
set of maximizer and minimizer (on the right of $\tmo$)
$\M_0\equiv\left\{\tm_0,\tM_r,\tm_r,\cdots,\tM_1,\tm_1,\tM_0\right\}$.

\noindent In the same way we construct the set $\M_0'$ by making
equivalent
 refinement  on the left of the valley $\{\tM_0',\tm_0,\tM_0\}$.
 We make a first refinement that gives the points $\{\tm_1',\tM'_1\}$,
 then we refine $\{\tM'_1,\tm_0\}$ and so on until we obtain
 $\{\tm_{r'}',\tM_{r'}'\}$ such that $\tm_0-\tM_{r'-1}' \geq b_n l_n$
and  $\tm_0-\tM_{r'}' \leq  b_n l_n$  (we denote $\M_0'$ this set
of maximizer and minimizer on the left of  $\tmo$). Finally we get
a set of maximizer and minimizer $\M\equiv \M_0'\cup \M_0=\{
\tM'_0
,\tm_1',\tM_1',\cdots,\tM'_{r'},\tm_0,\tM_{r},\cdots,\tM_1,\tm_1,\tM_0
\}$. \\

\noindent We will use the following notations,
\begin{eqnarray} \begin{array}{l|l}
 \textrm{If } 0 \leq i,j  \leq r  & \textrm{If } 0 \leq i,j \leq r'  \\
 \delta_{i,j}=S^{n}_{\tM_{i}}-S^{n}_{\tm_{j}}, & \delta_{i,j}'=S^{n}_{\tM_{i}'}-S^{n}_{\tm_{j}'}, \\
\eta_{i,j}=S^{n}_{\tM_{i}}-S^{n}_{\tM_{j}}, &
 \eta_{i,j}'=S^{n}_{\tM_{i}'}-S^{n}_{\tM_{j}'}, \\
\mu_{i,j}=S^{n}_{\tm_{i}}-S^{n}_{\tm_j}. &
\mu_{i,j}'=S^{n}_{\tm_{i}'}-S^{n}_{\tm_j'}. \\
\end{array}
\label{notsuper} \label{notprofvallée}
\end{eqnarray}
The beauty of the refinement is that we get immediately the
following relations between the random variables defined in
\ref{notsuper}
\begin{eqnarray}
& & \delta_{0,0} > \delta_{1,1} > \cdots > \delta_{r,r}\geq 0 ,
\\ & & \delta_{1,0}> \delta_{2,1} > \cdots > \delta_{r,0}\geq 0 ,
\end{eqnarray}
in the same way
\begin{eqnarray}
& & \delta'_{0,0} > \delta'_{1,1} > \cdots > \delta'_{r,r} \geq 0
,
\\ & & \delta'_{1,0}> \delta'_{2,0} > \cdots > \delta'_{r',0} \geq
0,
\end{eqnarray}
and
\begin{eqnarray}
\forall i,\ 0 \leq i \leq r-1,\  \eta_{i,i+1} \geq 0 , \\
\forall i,\ 0 \leq i \leq r'-1,\  \eta_{i,i+1}' \geq 0 .
\end{eqnarray}

\noindent We remark that the construction we made is possible if
and only if $\tm_0-\tM_{0}' \geq b_n l_n$ and $\tM_{0}-\tm_{0}
\geq l_n b_n $, but this is true with probability very near one,
indeed the following lemma will be proved in the Appendix
\ref{sec6} :

\begin{Lem} \label{interminib} There exists $h>0$ such that if \ref{hyp1bb}, \ref{hyp0} and \ref{hyp4} hold, for all $\gamma>0$ there exists $n_0 \equiv
n_0\left(\gamma \right)$ such that for all $n>n_0$, we have
\begin{eqnarray}
& & Q\left[\tM_0-\tmo \geq (\log n)^2 (65 \sigma^2 \log_2
n)^{-1}\right]\geq 1-h \left(( \log_3  n) (\log_2
n)^{-1}\right)^{1/2}
  \label{intermini1b} \\
& & Q\left[\tm_0-\tM_0' \geq (\log n)^2 (65 \sigma^2 \log_2
n)^{-1} \right]\geq 1- h \left(( \log_3  n) (\log_2
n)^{-1}\right)^{1/2}. \label{intermini2b}
\end{eqnarray}
\end{Lem}


\subsection{Definition of the set of good environments \label{sec4.5} }

\noindent Before defining a good environment, we introduce the
following random variables, let $\gamma>0$ and $n>3$,
\begin{eqnarray}
& & \tM_<=\sup\left\{m \in \Z,\ m < \tmo,\ S^n_{m}-S^n_{\tmo} \geq  \left(\log ( q_n   (\log n)^{\gamma})\right) (\log n )^{-1} \right\} , \label{Msup} \\
& &\tM_>=\inf\left\{m \in \Z,\ m> \tmo,\ S^n_{m}-S^n_{\tmo} \geq
\left(\log ( q_n   (\log n)^{\gamma})\right) (\log n )^{-1}
\right\} , \nonumber
\end{eqnarray}
where $q_n=\exp\left\{ \left( (200\sigma)^2 \gamma ( \log_2
n)^{7/2}(\log n)^{3/2}\right)^{1/2}\right\}$.

\begin{Rem} Proposition \ref{lem2} shows that for the scale of time $n$, Sinai's walk is trapped in the basic
valley $\{\tM_0',\tmo,\tM_0\}$. In the same way we will prove that
starting from $\tmo$ with a scale of time $q_n$, Sinai's walk is
trapped in
 the valley $\{\tM_<,\tmo,\tM_>\}$. This argument will be used in
 the proof of Theorem \ref{thSinai1}.
\end{Rem}

\noindent \\ Now we can define what we call a \textit{ good
environment }

\begin{Def} \label{super} Let $n>3$, $\kappa \in ]0,k_+[$, $\gamma>0$
 and $\omega \in \Omega_1$, we will say that $\alpha \equiv
\alpha(\omega)$ is a  \textit{ good environment } if
 the sequence $(\alpha_i,\ i \in \Z)\equiv(\alpha_i(\omega),\ i \in \Z)$
 satisfies the properties \ref{existence} to \ref{superdeltar}
\begin{eqnarray}
& \bullet & \textrm{ The valley }\{\tM_0',\tmo,\tM_0\} \textrm{  exists : } \label{existence} \\
& & 0 \in [\tM_0',\tM_0], \\
& & \delta_{0,0}\geq 1 +\gamma(n),\  \delta_{0,0}' \geq 1 +\gamma(n) , \label{3.18} \\
& & \textrm{ If } \tmo>0, S_{\tM_0'}-\max_{0 \leq m \leq \tmo }
\left(S^n_m\right)\geq \gamma(n) , \label{supermaxx}\\
& & \textrm{ if } \tmo<0, S_{\tM_0}-\max_{\tmo \leq m \leq 0 }
\left(S^n_m\right)\geq \gamma(n) . \label{supermaxx2} \\
 & \bullet & \max_{\tM_0' \leq l \leq \tM_0}
\left( (\alpha_l)^{-1}\right) \leq
(\log n)^{\frac{6}{\kappa}} \label{supermax} ,\\
& & \max_{\tM_0' \leq l \leq \tM_0} \left((\beta_l)^{-1}\right)
\leq
(\log n)^{\frac{6}{\kappa}} \label{supermax1} .\\
& \bullet & \tM_0 \leq (\sigma^{-1}\log n)^2 \log_2 n, -\tM_0'
\leq (\sigma^{-1}\log n)^2 \log_2 n . \label{superint}
\\
& \bullet & \tM_{<} \geq \tmo-L_n , \quad \tM_{>} \leq
\tmo+L_n . \label{superMsup} \\
& \bullet & r  \leq 2 (\log n)^{1/2}(\gamma \log_2 n)^{-1/2} \label{superr} , \\
& &   r'  \leq 2 (\log n)^{1/2}(\gamma \log_2 n)^{-1/2}
\label{superrpp} . \\
&  \bullet & \textrm{ For all } \ 0 \leq i \leq r-1 \nonumber \\
& & \eta_{i,i+1} \geq \gamma(n) \label{supereta} ,
\\ & & \delta_{i+1,i+1}\geq \gamma(n) \label{superdelta} , \\
& & \mu_{i+1,0}\geq \gamma(n) \label{supermu} . \\
& \bullet &  \textrm{ For all } 0 \leq i \leq r'-1 \nonumber  \\
& & \eta'_{i,i+1} \geq \gamma(n) \label{superetap} ,\\
& & \delta'_{i+1,i+1}\geq \gamma(n)  \label{superdeltap} ,
 \\
& & \mu'_{i+1,0}\geq \gamma(n)  \label{supermup} .
 \\
& \bullet &  \delta_{1,1} \leq 1-\gamma(n) ,\label{superdelta1}  \\
& &  \delta'_{1,1} \leq 1-\gamma(n)  .  \label{superdelta1p} \\
& \bullet & \delta_{r,r} \leq  (\log q_n) (\log
n)^{-1} , \label{superdeltarp} \\
& & \delta_{r',r'}' \leq   (\log q_n) (\log n)^{-1}.
\label{superdeltar}
\end{eqnarray}
where $L_n=\left( 8 \log[ (\log n)^{\gamma}q_n]
\sigma^{-1}\right)^2 \log_2 n$ and recalling that $q_n=\exp\left\{
\left( (200\sigma)^2 \gamma  ( \log_2 n)^{7/2}(\log
n)^{3/2}\right)^{1/2}\right\}$, $\delta_{.,.}$, $\delta_{.,.}'$,
$\eta_{.,.}$, $\eta_{.,.}'$,
 $\mu_{.,.}$ and $\mu_{.,.}'$ are given by \ref{notprofvallée} and
$\gamma(n)= ( \gamma \log_2 n) (\log n)^{-1}$.
\end{Def}

\noindent We define the set of \textit{good environments} $G_n$ as
\begin{eqnarray}
G_n=\left\{\omega \in \Omega_1,\ \alpha(\omega) \textrm{ is a
''good'' environment } \right\}.
\end{eqnarray}

\begin{Rem}
We remark that a good environment $\alpha$ is such that the
different random variables $\tM_0,\tM_0',\tmo,$ $r,r',
\delta_{.,.},\delta_{.,.}'$, $\mu_{.,.}$ and $\mu'_{.,.}$ that
depends on $\alpha$ satisfy some properties in relation to
deterministic parameters like $n$,
$\gamma$, $\sigma$ and $\kappa$. \\
\noindent The properties \ref{existence}-\ref{supermaxx2} concern
the existence of the basic valley $\{\tM_0',\tmo,\tM_0\}$ with his
main properties. \\
\noindent The properties \ref{supermax} and \ref{supermax1} are
technical properties due to the hypothesis \ref{hyp4}. There is no
equivalent properties in Sinai's paper because the
stronger hypothesis \ref{hyp2} is used. \\
\noindent \ref{superint} (respectively \ref{superMsup}) give an
upper bound of the distance between $\tM_0'$ and $\tM_0$
(respectively $\tM_<$ and $\tM_>$) and the origin (respectively to
the random point $\tmo$). \\
 \noindent   The properties from
\ref{superr} to \ref{superdeltar} concern the properties of the
valleys obtained by the ordered chopping
 of $\{\tM_0',\tmo,\tM_0\}$ effectuated in the previous paragraph.
 We remark that \ref{superr} and \ref{superrpp} give a deterministic upper
 bound for the number of right (respectively left) refinement
 performed in the ordered chopping in valleys, these upper bounds depend on
 $n$. This $n$ dependance that does not appear in Sinai's work comes from the fact that we perform a chopping in valleys
 in such a way that the successive valleys are nested and contain $\tmo$. This is a basic ingredient to get a result stronger than Sinai's one for the random walk itself.
\end{Rem}

\begin{Pro} \label{profonda}  There exists $h>0$ such that if \ref{hyp1}, \ref{hyp0} hold and for all $\kappa \in ]0,\kappa^+[$ \ref{hyp4} hold, for all $\gamma>0$,
 there exists $n_0 \equiv n_0(\kappa,\gamma) $ such that for all $n>n_0$
\begin{eqnarray}
Q\left[G_n\right] \geq 1-h \left((\log_3  n) (\log_2
n)^{-1}\right)^{1/2} .
\end{eqnarray}
\end{Pro}
\begin{Pre}
The proof of this proposition is done in the Appendix \ref{sec6}.
In fact $n_0 \equiv
n_0(\kappa,\gamma,\sigma,\E\left[|\epsilon_0|^3\right],\E\left[\epsilon_0^4\right],C)$,
 where $C=\E_Q\left[e^{\kappa \epsilon_0}\right] \vee
\E_Q\left[e^{-\kappa \epsilon_0}\right]$ but for simplicity we do
not always make explicit the dependance on $\sigma,
\kappa,\E\left[|\epsilon_0|^3\right],\E\left[\epsilon_0^4\right]$
and $C$ of $n_0$.
\end{Pre}



\section{Proof of the main results (Proposition \ref{lem2} and Theorem \ref{thSinai1}) \label{sec5}}

\subsection{Basic results for birth and death processes }

For completeness we recall some results of \cite{Chung} on
inhomogeneous discrete time birth and death processes, we will
always assume that $\alpha$ is fixed (denoted $\alpha \in
\Omega_1$ in this work).

\noindent Let $x$, $a$ and $b$ in $\Z$, $a\neq b$, suppose
$X_0=a$, denote

\begin{eqnarray}
& & T_b^a=\left\{\begin{array}{l}  \inf \{k\in\N^*,\ X_k=b \},  \\
 + \infty \textrm{, if such a } k \textrm{ not exists.}
\end{array} \right.
\end{eqnarray}
Assume $a<x<b$, the two following lemmata can be found in
\cite{Chung} (pages 73-76), their proof follow from the method of
difference equations.
\begin{Lem} \label{3.7bb} For all $\alpha \in \Omega_1$, we have
 \begin{eqnarray}
 & &
 \p^{\alpha}_x\left[T^x_a>T^x_b\right]=\frac{\sum_{i=a+1}^{x-1}\exp
 \left( \log n
\big(S_{i}^n-  S_{a}^n\big) \right)  +1}{\sum_{i=a+1}^{b-1}\exp
\left(\log n \big( S_{i}^n- S_{a}^n  \big) \right)+1 } \label{k1}
,
\\  & & \p^{\alpha}_x \left[T^{x}_a<T^{x}_b
\right]=\frac{\sum_{i=x+1}^{b-1} \exp \left( \log n \big( S_{i}^n-
S_{b}^n \big) \right) +1}{\sum_{i=a+1}^{b-1}\exp \left( \log n
\big( S_{i}^n- S_{b}^n \big) \right) +1 } \label{k2} .
\end{eqnarray}
\end{Lem}
\noindent Let us denote $T^x_a \wedge T^x_b $ the minimum between
$T^x_a$ and
 $T^x_b$.
\begin{Lem} \label{lA0} For all $\alpha \in \Omega_1$, we have
\begin{eqnarray}
&& \Ea_{a+1}\left[T^{a+1}_{a} \wedge T^{a+1}_{b} \right]=\frac{\sum_{l=a+1}^{b-1}\sum_{j=l}^{b-1}\frac{1}{\alpha_l}F_n(j,l)}{\sum_{j=a+1}^{b-1}F_n(j,a)+1}\label{c4s1l1e1} , \\
&& \Ea_x \left[ T^{x}_{a} \wedge
T^{x}_{b}\right]=\E^{\alpha}_{a+1}\left[T^{a+1}_{a} \wedge
T^{a+1}_{b}
\right]\left(1+\sum_{j=a+1}^{x-1}F_n(j,a)\right)-\sum_{l=a+1}^{x-1}\sum_{j=l}^{x-1}\frac{1}{\alpha_l}F_n(j,l),
\label{c4s1l1}
\end{eqnarray}
where  $F_n(j,l)=\exp \left( \log n \big( S_{j}^n- S_{l}^n \big)
\right) $.
\end{Lem}


\subsection{Proof of the sub-diffusive behavior (Proposition \ref{lem2} )}
\noindent



\noindent \textbf{Ideas of the proof} First we prove that starting
from $0$ the probability to hit $\tmo$ before one of the points
$\tM_0'-1$ or $\tM_0+1$ goes to $1$ (lemma \ref{lem1lem2}) and
starting from $\tm_0$ the probability
  of staying in
the interval $[\tM_0',\tM_0]$ in a time $n$ goes to $1$ when $n$
goes to
 infinity
 (lemma  \ref{lem2lem2}).
\\
\noindent
In this section we will always assume that $m_0<0$, (computations are the same for the other case). \\

\begin{Lem} \label{lem1lem2}  There exists $h>0$ such that if \ref{hyp1} and  \ref{hyp0} hold and for all $\kappa \in ]0,\kappa^+[$ \ref{hyp4} holds, for all
 $\gamma>2$ there exists $n_0\equiv
n_0(\gamma,\kappa)$ such that for all $n>n_0$ there exists $G_n
\subset \Omega_1$ with  $Q\left[G_n\right]\geq 1-h \left((\log_3
n) (\log_2 n)^{-1}\right)^{1/2} $ and for all $\alpha \in G_n$
\begin{eqnarray}
\p^{\alpha}_0\left[T^{\tilde{0}}_{\tmo} \geq
T^{\tilde{0}}_{\tM_0+1 }\right] & \leq & \sigma^{-2} (\log_2 n)
(\log n)^{-\gamma+2}+(n (\log n)^{\gamma})^{-1} . \label{Last1}
\end{eqnarray}
\end{Lem}
\begin{Pre}
Assume $\gamma>2$, using lemma \ref{3.7bb}  we easily get that
\begin{eqnarray*}
\p_0^{\alpha}\left[T^{\tilde{0}}_{\tmo} \geq T^{\tilde{0}}_{\tM_0
+1}\right]& \leq & |\tmo| \max_{\tmo+1 \leq i \leq
-1}\left(\exp\Big(-\log n \big(S^{n}_{\tM_0}-S^{n}_{i}\big) \Big)
\right)+1
\end{eqnarray*}
\noindent Using \ref{supermaxx2} and \ref{superint}, we get
\ref{Last1}
\end{Pre}


\noindent
\begin{Rem} \label{Remqp3}   By hypothesis $\tM_0'<\tm_0<0$ therefore $\p^{\alpha}\left[T^{\tilde{0}}_{\tmo}>T^{\tilde{0}}_{\tM_0'-1
}\right]=0$.
\end{Rem}

\begin{Lem} \label{lem2lem2} There exists $h>0$ such that if \ref{hyp1} and  \ref{hyp0} hold and for all $\kappa \in ]0,\kappa^+[$ \ref{hyp4} holds, for all
 $\gamma>2$ there exists $n_0\equiv
n_0(\gamma,\kappa)$ such that for all $n>n_0$ there exists $G_n
\subset \Omega_1$ with  $Q\left[G_n\right]\geq 1- h \left((\log_3
n) (\log_2 n)^{-1}\right)^{1/2} $ such that for all $\alpha \in
G_n$ we have
\begin{eqnarray}
& & \p^{\alpha}_{\tm_0}\left[T^{\tmo}_{\tM_0'-1} \wedge
T^{\tmo}_{\tM_0+1}>n \right] \geq  1- (\log n)^{-\gamma} ,
 \label{Last2}
 \end{eqnarray}
moreover
\begin{eqnarray}
& & \p^{\alpha}_{\tm_0}\left[T^{\tmo}_{-[ (\sigma^{-1} \log n)^2
\log_2 n ]-1} \wedge T^{\tmo}_{[ (\sigma^{-1}\log n)^2 \log_2
n]+1}>n \right] \geq  1- (\log n)^{-\gamma}  . \ \label{Last3}
 \end{eqnarray}
\end{Lem}
\begin{Pre}
\noindent For all $i\geq 2$, define
\begin{eqnarray}
&& T_i^{x\rightarrow x}= \left\{\begin{array}{l}  \inf\{k>T_{i-1},\ X_t=x\}, \\
 + \infty \textrm{, if such } k \textrm{ does not exist.}
 \end{array} \right. \\
& & T_1^{x\rightarrow x} \equiv T^{x\rightarrow
x}=\left\{\begin{array}{l}  \inf\{k\in \N^*,\ X_k=x
\textrm{ with } X_0=x \}, \\
 + \infty \textrm{, if such } k \textrm{ does not exist.}
 \end{array} \right.
\end{eqnarray}
We  denote $\tau_1=T_{1}^{x\rightarrow x}$ and
$\tau_i=T_{i}^{x\rightarrow x}-T_{i-1}^{x\rightarrow x}$, for all
$i\geq 2$. Let $n\geq 1$, remark that $T_{n}^{\tmo \rightarrow
\tmo} \equiv \sum_{i=1}^{n}\tau^{\tmo \rightarrow \tmo}_i
 > n$ so
\begin{eqnarray}
\p^{\alpha}_{\tm_0}\left[ T^{\tmo}_{\tM_0'-1} \wedge
T^{\tmo}_{\tM_0+1}>n \right]& =&
\p^{\alpha}_{\tm_0}\left[T^{\tmo}_{\tM_0'-1} \wedge
T^{\tmo}_{\tM_0+1}>n,\sum_{i=1}^{n}\tau^{\tmo \rightarrow \tmo}_i
> n  \right] \\
&\geq & \p^{\alpha}_{\tm_0}\left[T^{\tmo}_{\tM_0'-1} \wedge
T^{\tmo}_{\tM_0+1} > \sum_{i=1}^{n}\tau^{\tmo \rightarrow \tmo}_i
\right]  \label{probi}
\end{eqnarray}
By the strong Markov property the random variables $(\tau_i,1 \leq
i\leq n)$ are i.i.d therefore
\begin{eqnarray}
\p^{\alpha}_{\tm_0}\left[T^{\tmo}_{\tM_0'-1} \wedge
T^{\tmo}_{\tM_0+1} > \sum_{i=1}^{n}\tau^{\tmo \rightarrow \tmo}_i
\right] = \left(\p^{\alpha} \left[ T^{\tmo \rightarrow \tmo} \leq
T^{\tmo}_{\tM_0'-1}\wedge T^{\tmo}_{\tM_0+1}\right]\right)^{n} .
\label{lc1}
\end{eqnarray}
Moreover it is easy to check that
\begin{eqnarray}
& \qquad &  \p^{\alpha}_{\tm_0} \left[ T^{\tmo \rightarrow \tmo}
\leq  T^{\tmo}_{\tM_0'-1}\wedge
T^{\tmo}_{\tM_0+1}\right]=\alpha_{\tmo}\p^{\alpha}_{\tmo+1} \left[
T^{\tmo+1}_{\tM_0+1} \leq
T^{\tmo+1}_{\tmo}\right]+\beta_{\tmo}\p^{\alpha}_{\tmo-1}
\left[T^{\tmo-1}_{\tM_0'-1} \leq T^{\tmo-1}_{\tmo} \right] .
\label{eqquipu1}
\end{eqnarray}
 \noindent Using \ref{k1} and \ref{3.18}  we get that there exists $n_0 \equiv
n_0(\kappa,\gamma)$ such that for all $n>n_0$ and  all $ \alpha
\in G_n$, $ \p^{\alpha}_{\tm_0+1}\left[T^{\tmo +1}_{\tM_0+1}<
\right.$ $ \left. T^{\tmo+1}_{\tmo}\right]  \leq
n^{-(1+\gamma(n))} $, in the same way $
\p^{\alpha}_{\tm_0-1}\left[T^{\tmo-1}_{\tM_0'-1}<T^{\tmo-1}_{\tmo}\right]
  \leq n^{-(1+\gamma(n))} $.
Using this and \ref{eqquipu1}, we get for  $n>n_0$ and  all
$\alpha \in G_n$
\begin{eqnarray}
& &\p^{\alpha}_{\tm_0}\left[T^{\tmo}_{\tM_0'-1}\wedge
T^{\tmo}_{\tM_0+1}<T^{ \tmo\rightarrow \tmo}\right]\leq
n^{-1-\gamma(n)} . \label{theq1213}
\end{eqnarray}
\noindent Replacing \ref{theq1213} in \ref{lc1} and  using
\ref{probi} and the fact $\left(1-x \right)^n \geq 1-n x$, for all
$0 \leq x \leq
 1 $ and  all $n\geq 1$ we get \ref{Last2}. For \ref{Last3} we
 use \ref{Last2} and \ref{superint}.
\end{Pre}
\\
\begin{Prepr}{\ref{lem2}}
\noindent By the strong Markov property and remark \ref{Remqp3} we
get that
\begin{eqnarray}
& &
\p^{\alpha}_0\left[\bigcap_{k=0}^n\left\{X_m\in\left[\tM_0',\tM_0
\right]\right\}\right] \geq
\p^{\alpha}_{\tm_0}\left[T^{\tmo}_{\tM_0'-1} \wedge T^{\tmo
}_{\tM_0+1 } >
n\right]-\p_0^{\alpha}\left[T^{\tilde{0}}_{\tmo}>T^{\tilde{0}}_{\tM_0+1
}\right],
\end{eqnarray}
Using Lemmata \ref{lem1lem2} and  \ref{lem2lem2}, we get
\ref{2.19}. We get \ref{2.20} using \ref{2.19} and \ref{superint}.
\end{Prepr}

\noindent \\ The next lemma will be used for the proof of Theorem
\ref{thSinai1}.

\begin{Lem} \label{lemsauve} There exists $h>0$, such that if \ref{hyp1} and  \ref{hyp0} hold and for all $\kappa \in ]0,\kappa^+[$ \ref{hyp4} holds, for all
 $\gamma>2$ there exists $n_0\equiv
n_0(\gamma,\kappa)$ such that for all $n>n_0$ there exists $G_n
\subset \Omega_1$ with  $Q\left[G_n\right]\geq 1- h \left((\log_3
n) (\log_2 n)^{-1}\right)^{1/2} $ and for all $\alpha \in G_n$ we
have
\begin{eqnarray}
& &\p^{\alpha}_{\tm_0}\left[T^{\tmo}_{\tmo-L_n} \wedge T^{\tmo
}_{\tmo+ L_n}>q_n \right] \geq  1 -(\log n)^{-\gamma} \ ,
 \end{eqnarray}
where $L_n$ and $q_n$ are given
 at the end of Definition \ref{super}.
\end{Lem}

\begin{Pre}
Using what we did to prove Lemma \ref{lem2lem2} replacing $\tM_0$
by $\tM_>$ and $\tM_0'$ by $\tM_<$ (see \ref{Msup} for the
definitions of $\tM_>$ and $\tM_>$), we easily get this lemma.
\end{Pre}

\subsection{Analysis of the return time $T^{\tmo \rightarrow
\tmo}$ \label{3par2}}

It is easy to check that  $\Eam\left[T^{\tmo \rightarrow
\tmo}\right]=\infty$ $Q.a.s$,  however we will need an upper bound
for the probability $\pam\left[T^{\tmo \rightarrow \tmo}>k\right]$
with $k>0$. We denote $a \vee b = \max(a,b)$.

\begin{Lem} \label{lemsecmom1} For all $\alpha \in \Omega_1$ and  all $n>1$, we
have for all $i,\ 0 \leq i \leq r$
\begin{eqnarray}
&&
\Ea_{\tmo+1}\left[\left(T^{\tilde{m}_{0}+1}_{\tilde{m}_{0}}\wedge
T^{\tilde{m}_{0}+1}_{\tM_i+1}\right)^2 \right] \leq D_i
n^{(\delta_{i+1,i+1}-\eta_{i,i+1})\vee 0}, \label{eqlemsecmom1}
\end{eqnarray}
 with $D_i\equiv D_i(\alpha,n)=|\tM_{i}-\tm_0|^5 \left(\max_{ \tmo \leq
l \leq \tM_{i}} \left(\frac{1}{\alpha_l}\right) \right)^2$, and
for all $i,\ 0 \leq i \leq r'$
\begin{eqnarray}
\Ea_{\tmo-1}\left[\left(T^{\tilde{m}_{0}-1}_{\tilde{m}_{0}}\wedge
T^{\tilde{m}_{0}-1}_{\tM_{i}'-1}\right) ^2 \right] \leq D_i'
 n^{(\delta_{i+1,i+1}'-\eta_{i,i+1}')\vee 0}, \label{eqlemsecmom2}
\end{eqnarray}
with $D_i'\equiv D_i'(\alpha,n)=|\tM_{i}'-\tm_0|^5 \left(\max_{
\tM_{i}' \leq l \leq \tmo} \left(\frac{1}{\beta_l}\right)
\right)^2$. See \ref{notprofvallée} for the definitions of
$\eta_{i,i+1}',\ \delta_{i+1,i+1}',\ \eta_{i,i+1}$ and
 $\delta_{i+1,i+1}$, recalling that $r$ and $r'$ are (respectively) the number of right (respectively left)
 refinement (see section \ref{OC}).
\end{Lem}
\begin{Pre}
 We only prove \ref{eqlemsecmom1} ( the proof of \ref{eqlemsecmom2} is identical). It is easy to check, with the method of difference equations,
\begin{eqnarray}
\Ea_{\tmo}\left[\left(T^{\tmo}_{\tmo+1}\wedge
T^{\tmo}_{\tM_{i}+1}\right)^2\right]=\frac{\sum_{l=\tmo+1}^{\tM_{i}}\sum_{j=\tmo+1}^l\frac{2u_l-1}{\alpha_l}F_n(j,l)}{\sum_{j=\tmo+1}^{\tM_{i}}F_n(j,\tmo)+1}
, \label{eqsecmom2}
\end{eqnarray}
with
\begin{eqnarray}
&& u_l  = \Ea_l\left[T^{l}_{\tmo} \wedge T^l_{\tM_i+1} \right],
\label{eqsecmom48}
\end{eqnarray}
$u_l$ is given by \ref{c4s1l1} and $F_n(.,.)$ at the end of  Lemma
\ref{lA0}. First we give an upper bound of \ref{eqsecmom48}.
Denoting $C_i \equiv C_i(\alpha,n)= \max_{ \tmo \leq l \leq
\tM_{i}} \left( \frac{1}{ \alpha_l }\right)(\tM_{i}-\tm_0)^2$ it
is easy to check that $u_l \leq C_i
\left(1+\sum_{j=\tmo+1}^{l-1}F_n(j,\tmo)\right) $. We have
\begin{eqnarray}
& &
\sum_{l=\tmo+1}^{\tM_{i}}\sum_{j=\tmo+1}^l\frac{2u_l-1}{\alpha_l}F_n(j,l)
  \leq   2C_i
\sum_{l=\tmo+1}^{\tM_{i}}\sum_{j=\tmo+1}^l
 \left(1+\sum_{i=\tmo+1}^{l-1}F_n(i,\tmo)\right) (\alpha_l)^{-1}F_n(j,l)
 .
 \label{seceq1}
\end{eqnarray}
Now let us consider the first refinement of $\{\tmo,\tM_{i}\}$,
 denote $\tm_{i+1}$ the minimizer obtained and $\tM_{i+1}$
the  maximizer, it is easy to check (see Figure \ref{thfig3}) that
\begin{eqnarray}
 & & \sum_{l=\tmo+1}^{\tM_{i}}\sum_{j=\tmo+1}^l
\frac{\left(1+\sum_{i=\tmo+1}^{l-1}F_n(i,\tmo)\right)}{\alpha_l}F_n(j,l)
  \leq
 \frac{|\tM_{i}-\tm_0|^3}{2} \max_{ \tmo \leq l \leq
\tM_{i}} \left(\frac{1}{\alpha_l}\right)  n^{(\delta_{i,0}) \vee
(\delta_{i+1,0}+\delta_{i+1,i+1}) } \label{seceq2} ,
\end{eqnarray}
where $\delta_{.,.}$ is given in \ref{notsuper}. Using
\ref{seceq1} and \ref{seceq2} we get
\begin{eqnarray}
&&
\sum_{l=\tmo+1}^{\tM_{i}}\sum_{j=\tmo+1}^l\frac{2u_l-1}{\alpha_l}F_n(j,l)
\leq  D_i\times n^{(\delta_{i,0}) \vee
(\delta_{i+1,0}+\delta_{i+1,i+1}) } ,\label{eqsecfin1}
\end{eqnarray}
\noindent where $D_i\equiv D_i(\alpha,n)=|\tM_{i}-\tm_0|^5
\left(\max_{ \tmo \leq l \leq \tM_{i}}
\left(\frac{1}{\alpha_l}\right) \right)^2$.


\noindent Moreover it is easy to check that $
\sum_{j=\tmo+1}^{\tM_{i}}F_n(j,\tmo) \geq n^{\delta_{i,0}}$,
replacing this and \ref{eqsecfin1}  in \ref{eqsecmom2} and
noticing that $\delta_{i+1,0}-\delta_{i,0}=-\eta_{i,i+1}$
 we get \ref{eqlemsecmom1}.
\end{Pre}

\begin{Pro} \label{lemtrm1} For all $\alpha \in \Omega_1$, $n>1$ and $q>0$ we
have, for all $i$, $0 \leq i \leq r$
\begin{eqnarray}
\p^{\alpha}_{\tmo+1}\left[T^{\tmo+1}_{ \tmo }
> q \right] \leq (D_i n^{(\delta_{i+1,i+1}-\eta_{i,i+1})\vee 0})q^{-2} +
n^{-\delta_{i,0}} ,    \label{eqlemtrm2}
\end{eqnarray}
with $D_i=|\tM_{i}-\tm_0|^5 \left(\max_{ \tmo \leq l \leq
\tM_{i}} \left(\frac{1}{\alpha_l}\right)\right)^2$, and for all
$i$, $ 0 \leq i \leq r'$
\begin{eqnarray}
\p^{\alpha}_{\tmo-1}\left[T^{\tmo-1}_{\tmo} > q \right] \leq (D_i'
n^{(\delta_{i+1,i+1}'-\eta_{i,i+1}')\vee 0})q^{-2} +
n^{-\delta_{i,0}'} \ \label{eqlemtrm3} .
\end{eqnarray}
with $D_i'=|\tM_{i}'-\tm_0|^5 \left(\max_{ \tM_{i}' \leq l \leq
\tmo} \left(\frac{1}{\beta_l}\right) \right)^2 $. See
\ref{notprofvallée} for the definitions of  $\eta_{i,i+1}',\
\delta_{i+1,i+1}',\ \eta_{i,i+1}$ and
 $\delta_{i+1,i+1}$, recalling that $r$ and $r'$ are (respectively) the number of right (respectively left)
 refinement (see section \ref{OC}).
\end{Pro}

\begin{Rem}
\ref{eqlemtrm2} does not imply that
$\p^{\alpha}_{\tmo+1}\left[T^{\tmo+1}_{ \tmo }
> q \right]$ is sumable on $q$, indeed on the right
hand side of \ref{eqlemtrm2}, "$n^{-\delta_{i,0}}$" does not
depend on $q$.
\end{Rem}

\begin{Prepr}{\ref{lemtrm1}}
Let us estimate
$\p^{\alpha}_{\tmo+1}\left[T^{\tilde{m}_{0}+1}_{\tilde{m}_{0}}
> q\right]$, let $0 \leq i \leq r$,  we have
\begin{eqnarray}
\pa_{\tmo}\left[T^{\tilde{m}_{0}+1}_{\tilde{m}_{0}}>q\right] \leq
\pa_{\tmo+1}\left[T^{\tilde{m}_{0}+1}_{\tilde{m}_{0}} \wedge
T^{\tilde{m}_{0}+1}_{\tM_i+1}>q\right]  +
\p^{\alpha}_{\tmo+1}\left[
T^{\tilde{m}_{0}+1}_{\tilde{m}_{0}}>T^{\tilde{m}_{0}+1}_{\tM_i+1}\right]
. \label{eqpprobmom}
\end{eqnarray}
Using \ref{k2} and recalling that
$\delta_{i,0}=S^{n}_{\tM_i}-S^{n}_{\tmo}$ we get $
\p^{\alpha}_{\tmo+1}\left[
T^{\tilde{m}_{0}+1}_{\tilde{m}_{0}}>T^{\tilde{m}_{0}+1}_{\tM_i+1}\right]\leq
n^{-\delta_{i,0}} $. Moreover, by Markov inequality we have
$\p^{\alpha}_{\tmo+1}\left[T^{\tilde{m}_{0}+1}_{\tilde{m}_{0}}
\wedge T^{\tilde{m}_{0}+1}_{\tM_i+1}>q\right]\leq \left
(\Ea_{\tmo+1} \left[\left(T^{\tilde{m}_{0}+1}_{\tilde{m}_{0}}
\wedge T^{\tilde{m}_{0}+1}_{\tM_i+1} \right)^2\right] \right)
q^{-2}$ To end the proof we use \ref{eqlemsecmom1} (similar
computations give \ref{eqlemtrm3}).
\end{Prepr}

\subsection{Proof of Theorem \ref{thSinai1} \label{sec5.2}}


The sketch of the proof is the following we prove (with a
probability very near one) that $(X_k)_{1 \leq k \leq n}$ hit
$\tmo$ in a time smaller than $ n$. Then we show that it does not
exist an
 instant $1 \leq k \leq n-q_n$ ($q(n)$ is given at the end of Definition \ref{super}) such that the
\rw will not return to $\tmo$ (Proposition \ref{ProS1}). Finally
we prove that starting from $\tm_0$, in a time smaller than
$n-(n-q_n)=q_n$ the \rw can not escape from a region which size is
 of order $(\log q_n)^2$ (Proposition \ref{ProS2}) .

\noindent \\ First we introduce the next event, let $n>1$ and  $1
\leq q \leq n$
\begin{eqnarray}
\A_q=\bigcup_{ n-q \leq k \leq n}\left\{X_k=\tmo \right\} .
\end{eqnarray}
Let $\delta_q>0$, we have
\begin{eqnarray}
\p^{\alpha}_0\left[\left|\frac{X_n}{(\log
n)^2}-m_{0}\right|>\delta_q\right] &\leq&
\p_0^{\alpha}\left[\left|\frac{X_n}{(\log
n)^2}-m_{0}\right|>\delta_q,\
\A_q\right]+\p^{\alpha}_0\left[\A_q^c\right]. \label{ladereqs}
\end{eqnarray}
Now we estimate each probability of the right hand side of
\ref{ladereqs} in Propositions \ref{ProS1} and \ref{ProS2}.

\begin{Pro} \label{ProS1} There exists $h>0$ such that if \ref{hyp1} and  \ref{hyp0} hold
and for all $\kappa \in ]0,\kappa^+[$ \ref{hyp4} holds, for all
 $\gamma>12/\kappa+21/2$ there exists $n_0\equiv
n_0(\gamma,\kappa)$ such that for all $n>n_0$ there exists $G_n
\subset \Omega_1$ with  $Q\left[G_n\right]\geq 1-h \left( ( \log_3
n)/( \log_2 n)\right)^{1/2} $ and for all $\alpha \in G_n$
\begin{eqnarray}
\p^{\alpha}_0\left[\A_{q_n}^c\right] \leq \frac{2(\log_2
n)^{9/2}}{(\gamma )^{1/2}(\log
n)^{\gamma-(12/\kappa+21/2)}}+\Er\left(\frac{ (\log_2)^2}{(\log
n)^{\gamma-(6/\kappa+4)}} \right), \label{eqProS11}
\end{eqnarray}
$q_n$ is given at the end of Definition \ref{super}.
\end{Pro}

\noindent
\begin{Pre}
First we remark that for all $n>1$ and  all $ 1 \leq q \leq n $
\begin{eqnarray}
\p^{\alpha}_0\left[\A_{q}^c\right] & \leq & \p^{\alpha}_0\left[
T^{0}_{\tilde{m}_0}>n \right] + \p^{\alpha}_0\left[\A_{q}^c,\
T^{0}_{\tilde{m}_0}\leq n \right] \ . \label{6eq8}
\end{eqnarray}
We estimate each term of the right hand side of \ref{6eq8}, the
first one in Lemma \ref{lempro1} and the second in Lemma
\ref{lempro2}
\begin{Lem} \label{lempro1}There exists $h>0$ such that if \ref{hyp1} and  \ref{hyp0} hold
and for all $\kappa \in ]0,\kappa^+[$ \ref{hyp4} holds, for all
$\gamma > \frac{6}{\kappa}+4$, there exists $n_1' \equiv
n_1'(\kappa,\gamma)$ such that for all $n>n_1'$ there exists $G_n
\subset \Omega_1$ with $Q\left[G_n\right]\geq 1-h \left( ( \log_3
n)/( \log_2 n)\right)^{1/2} $ and for all $\alpha \in G_n$, we
have
\begin{eqnarray}
 \p^{\alpha}_0\left[ T^{0}_{\tilde{m}_{0}}>n \right]\leq
 \frac{5 (\log_2 n)^2}{\sigma^{4}(\log
 n)^{\gamma-\left(\frac{6}{\kappa}+4\right)}}.\label{4.36}
\end{eqnarray}
\end{Lem}
\begin{Pre}
Let us consider the valley $\{\tM_0',\tmo,\tM_0\}$, we assume
$\tm_0>0$ (computations are similar if $\tm_0 \leq 0$). We have
\begin{eqnarray}
 \p^{\alpha}_0\left[T^{0}_{\tmo}>n\right] &\leq& \p_0^{\alpha}\left[T^{0}_{\tmo} \wedge
T^{0}_{\tM_{0}'-1}>n \right] +
\p^{\alpha}_0\left[T^{0}_{\tM_{0}'-1}<T^{0}_{\tmo} \right] \
.\label{lasom1b}
\end{eqnarray}
For the second probability on the right hand side of \ref{lasom1b}
we have already see (lemme \ref{lem1lem2}) that for
 all $\gamma>2$ there exists $n_1\equiv
n_1(\kappa,\gamma)$ such that for all $n>n_1$ and all $\alpha \in
G_n$
\begin{eqnarray}
 \p^{\alpha}_0\left[T^{0}_{\tM_{0}'-1}<T^{0}_{\tmo} \right]\leq
 \sigma^{-2} \log_2 n (\log n)^{-\gamma+2} .\label{eqlempro1}
\end{eqnarray}
For the first probability on the right hand side of \ref{lasom1b}
we have by the Markov inequality
\begin{eqnarray}
\p_0^{\alpha}\left[T^{0}_{\tmo} \wedge T^{0}_{\tM_{0}'-1}>n
\right] \leq \E_0\left[T^{0}_{\tmo} \wedge
T^{0}_{\tM_{0}'-1}\right] n^{-1} .\label{4eq14}
\end{eqnarray}
To compute the mean in \ref{4eq14} we use lemma \ref{c4s1l1}, it
is easy to check that :
\begin{eqnarray}
 \E^{\alpha}_0\left[T^{0}_{\tM_{0}'-1} \wedge
T^{0}_{\tmo} \right] & \leq &
\sum_{l=\tM_{0}'}^{\tmo-1}\sum_{j=l}^{\tmo-1}\frac{1}{\alpha_l}F_n(j,l)
 \label{eqleprod1}
\end{eqnarray}
where
 $F_n(j,l)=\exp \left(\log n(S^n_l-S^{n}_j)\right)$. Let us consider the first refinement of
$\{\tM_{0}',\tmo\}$, it gives the point  $\tM_{1}'$ (for the
maximizer) and  $\tilde{m}_1'$ (for the minimizer), so we get
\begin{eqnarray}
\sum_{l=\tM_{0}'}^{\tmo-1}\sum_{j=l}^{\tmo-1}\frac{1}{\alpha_l}F_n(j,l)
& \leq &  C_0  n^{\delta_{1,1}'} , \label{eqleprod3}
\end{eqnarray}
where  $\delta_{1,1}'\equiv S^{n}_{\tM_1'}-S^n_{\tm_1'}$ and $C_0
\equiv C_0(\alpha,n)=(\tM_{0}'-\tm_0)^2\max_{\tM'_0 \leq l \leq
\tmo }\left(\frac{1}{\alpha_l}\right)$. Using \ref{eqleprod3},
 \ref{eqleprod1} and \ref{4eq14} we get
\begin{eqnarray}
\p^{\alpha}_0\left[T^{0}_{\tmo} \wedge T^{0}_{\tM_{0}'-1}>n
\right] \leq (C_0 n^{\delta_{1,1}'}) n^{-1}  .
\end{eqnarray}
Using formulas \ref{supermax}, \ref{superint} and
\ref{superdelta1p} we get that for all $\gamma>
\frac{6}{\kappa}+4$, there exists $n_2 \equiv n_2(\gamma)$ such
that for all $n>n_2$ and $\alpha \in G_n$
\begin{eqnarray}
\p^{\alpha}_0\left[T^{0}_{\tmo} \wedge T^{0}_{\tM_{0}'}>n \right]
& \leq & \frac{(2 \log_2 n)^2}{\sigma^4 (\log
n)^{\gamma-\left(\frac{6}{\kappa}+4\right)}}. \label{eqlempro2}
\end{eqnarray}
 \noindent We get \ref{4.36} using \ref{lasom1b}, \ref{eqlempro1} and \ref{eqlempro2}
and  taking $n_1'=n_1 \vee n_2$.
\end{Pre}

\begin{Lem} \label{lempro2} There exists $h>0$, such that if \ref{hyp1} and  \ref{hyp0} hold
and for all $\kappa \in ]0,\kappa^+[$ \ref{hyp4} holds, for all
 $\gamma>12/\kappa+21/2$ there exists $n_0\equiv
n_0(\gamma,\kappa)$ such that for all $n>n_0$ there exists $G_n
\subset \Omega_1$ with  $Q\left[G_n\right]\geq 1-h \left( (\log_3
n) (\log_2 n)^{-1}\right)^{1/2} $ and for all $\alpha \in G_n$
\begin{eqnarray}
 \p^{\alpha}_0 \left[\A_{q_n}^c,\ T^{0}_{\tilde{m}_{0}}\leq n \right] \leq \frac{3(\log_2 n)^{9/2}}{\sigma^{10}(\gamma)^{1/2}(\log
n)^{\gamma-(\frac{12}{\kappa}+\frac{21}{2})}}+\Er\left(\frac{1}{(\log
n)^{\gamma-1/2}(\log_2 n)^{1/2}}\right) \label{eqlempro21}
 \end{eqnarray}
$q_n$ is given at the end of definition \ref{super}.
\end{Lem}
\begin{Pre}
We recall that for all $ 1\leq q \leq n$ we have denoted
$\A_q^c=\bigcap_{n-q \leq k \leq n}\left\{X_k \neq \tmo \right\}$.
Denoting
\begin{eqnarray}
&& \bar{\A}_q^c = \bigcup_{1 \leq p \leq n-q-1} \left\{
\left\{X_p=\tmo \right\}\bigcap_{m=p+1}^n \left\{\ X_m \neq \tmo
\right\}\right\}, \label{4eq32}
\end{eqnarray}
we remark that $\left\{\A_q^c, T^0_{\tmo} \leq n\right\} \subset
\bar{\A}_q^c $. Therefore we only have to give an upper bound of $
\p^{\alpha}_0 \left[\bar{\A}_q^c \right]$, by the Markov property
we have
\begin{eqnarray}
\p^{\alpha}_0\left[\bar{\A}_q^c \right] &=&  \sum_{1 \leq p \leq
n-q-1} \p^{\alpha}_{\tmo}\left[ \bigcap_{m= 1}^{n-p} \left\{\ X_m
\neq \tmo \right\} \right]\p^{\alpha}_0\left[X_p=\tmo \right] .
\end{eqnarray}
Using the change $k=n-p$, we get
\begin{eqnarray}
 \p^{\alpha}_0\left[\bar{\A}_q^c \right] & \leq & \sum_{q+1 \leq k \leq n-1}
\p^{\alpha}_{\tmo}\left[ \bigcap_{m= 1}^{k} \left\{\ X_m \neq \tmo
\right\} \right]  \equiv  \sum_{q+1 \leq k \leq n-1} \pam \left[
T^{\tilde{m}_0\rightarrow \tilde{m}_0}
> k  \right] . \label{thshypsecmom1}
\end{eqnarray}

\begin{Rem}  We recall that \rw is null recurrent
$\p.a.s$, so for the moment, we can't say anything on  $\sum_{q+1
\leq k \leq n-1}$ $\p^{\alpha}_{\tmo}\left[
T^{\tilde{m}_0\rightarrow \tilde{m}_0} > k \right]$.
\end{Rem}
\noindent First let us decompose the sum in \ref{thshypsecmom1}
\begin{eqnarray}
 \sum_{q+1 \leq k \leq n-1} \p^{\alpha}_{\tmo}\left[
T^{\tilde{m}_0\rightarrow \tilde{m}_0} > k \right]
 &=&  \sum_{q \leq k \leq n-2} \alpha_{\tmo}
\p^{\alpha}_{\tmo+1}\left[  T^{\tmo+1}_{\tmo} > k
 \right] \label{eqtamere1}\\
&+& \sum_{q \leq k \leq n-2}
\beta_{\tmo}\p^{\alpha}_{\tmo-1}\left[
T^{\tilde{m}_0-1}_{\tilde{m}_0}
> k  \right]. \label{eqtamere2}
\end{eqnarray}
Let us  give an upper bound to the sum on the right hand side of
\ref{eqtamere1}. We want to find $q$ as small as possible but such
that this sum goes to 0. For this we use step by step the
inequality \ref{eqlemtrm2} to
$\p^{\alpha}_{\tmo+1}\left[T^{\tilde{m}_0+1}_{\tilde{m}_0}
> k \right]$ :  we have
\begin{eqnarray}
 \sum_{ [n^{\delta_{r,r}}]+1 \leq k \leq n-2} \p^{\alpha}_{\tmo+1}\left[
T^{\tilde{m}_0+1}_{\tilde{m}_0}
> k
\right] &=&
\sum_{k=[n^{\delta_{1,1}}]+1}^{n-2}\p^{\alpha}_{\tmo+1}\left[
T^{\tilde{m}_0+1}_{\tilde{m}_0} > k  \right] \label{lasom1} \\
&+&
\sum_{i=1}^{r-1}\sum_{k=[n^{\delta_{i+1,i+1}}]+1}^{[n^{\delta_{i,i}}]}\p^{\alpha}_{\tmo+1}\left[
 T^{\tilde{m}_0+1}_{\tilde{m}_0} > k
\right] . \label{lasom2}
\end{eqnarray}
For the sum on the right hand side of \ref{lasom1}, by inequality
\ref{eqlemtrm2} (taking $i=0$) we have
\begin{eqnarray}
 \sum_{k=[n^{\delta_{1,1}}]+1}^{n-2} \p^{\alpha}_{\tmo}\left[
T^{\tilde{m}_0+1}_{\tilde{m}_0} > k  \right]  & \leq &
\frac{n-n^{\delta_{1,1}}}{n^{\delta_{0,0}}}+\sum_{k=[n^{\delta_{1,1}}]+1}^n\frac{D_{0}
n^{(\delta_{1,1}-\eta_{0,1})\vee 0}}{k^2}\\ & \leq &
\frac{n}{n^{\delta_{0,0}}}+\frac{D_{0} }{n^{\delta_{1,1}\wedge
\eta_{0,1}}} , \label{eqeqlempro3}
\end{eqnarray}
where $D_0=|\tM_0-\tm_0|^5 \left( \max_{\tm_0\leq l \leq \tM_0}
\left(\frac{1}{ \alpha_l}\right) \right)^2$. For the other terms
$(1 \leq i \leq r-1)$ of the sum in \ref{lasom2}, using the
inequality \ref{eqlemtrm2} we have
\begin{eqnarray}
 \sum_{k=[n^{\delta_{i+1,i+1}}]+1}^{[n^{\delta_{i,i}}]}
\p^{\alpha}_{\tmo+1}\left[ T^{\tilde{m}_0+1}_{\tilde{m}_0} > k
\right]  & \leq &
\frac{n^{\delta_{i,i}}-n^{\delta_{i+1,i+1}}}{n^{\delta_{i,0}}}+\sum_{k=[n^{\delta_{i+1,i+1}}]+1}^{[n^{\delta_{i,i}}]}\frac{D_{i}(n^{(\delta_{i+1,i+1}-\eta_{i,i+1})\vee 0}}{k^2}\\
&\leq & \frac{1}{n^{\mu_{i,0}}}+\frac{D_{i}
}{n^{\delta_{i+1,i+1}\wedge \eta_{i,i+1}}} , \label{theq2343}
\end{eqnarray}
where we have used that $\delta_{i,0}-\delta_{i,i}=\mu_{i,0}$
 and $D_i=|\tM_i-\tm_0|^5 \left( \max_{\tm_0\leq l \leq \tM_i}
\left(\frac{1}{ \alpha_l}\right) \right)^2$. So, for the sum
\ref{lasom2} we get from \ref{theq2343}
 that
\begin{eqnarray}
\sum_{i=1}^{r-1}
\sum_{k=[n^{\delta_{i+1,i+1}}]+1}^{[n^{\delta_{i,i}}]}
\p^{\alpha}_{\tmo+1}\left[ T^{\tilde{m}_0+1}_{\tilde{m}_0} > k
\right] &\leq &
\sum_{i=1}^{r-1}\frac{1}{n^{\mu_{i,0}}}+\sum_{i=1}^{r-1}\frac{D_{i}
}{n^{\delta_{i+1,i+1}\wedge \eta_{i,i+1}}}
\\
& \leq &  \frac{r-1}{n^{\min_{1 \leq i \leq r-1}
\left(\mu_{i,0}\right)}} + \frac{(r-1) D_0}{n^{\min_{1 \leq i \leq
r-1} \left(\delta_{i+1,i+1}\wedge \eta_{i,i+1} \right)}} ,
 \label{eqeqlempro2}
\end{eqnarray}
and we have used that $D_{i}$ is decreasing in $i$. Collecting the
terms \ref{eqeqlempro2} and  \ref{eqeqlempro3} we get
\begin{eqnarray}
 \sum_{[n^{\delta_{r,r}}]+1 \leq k \leq n-2} \alpha_{\tmo}
\p^{\alpha}_{\tmo+1}\left[  T^{\tilde{m}_0+1}_{\tilde{m}_0}
> k
\right] & \leq & \frac{n}{n^{\delta_{0,0}}}+\frac{r-1}{n^{\min_{1
\leq i \leq r-1} \left(\mu_{i,0}\right)}}  + \frac{r
D_0}{n^{\min_{0 \leq i \leq r-1} \left(\delta_{i+1,i+1}\wedge
\eta_{i,i+1} \right)}} . \label{eqdelafin}
\end{eqnarray}
Now using the good properties \ref{3.18}, \ref{supermax},
\ref{supereta}, \ref{superdelta}, \ref{supermu}, \ref{superint}
and \ref{superr} we easily get that for all $\gamma
> \frac{12}{\kappa}+\frac{21}{2}$, there exist $n_1$ such that for all $n>n_1$, $\alpha \in G_n$,
\begin{eqnarray}
  \sum_{[n^{\delta_{r,r}}]+1 \leq k \leq n-2} \alpha_{\tmo}
\p^{\alpha}_{\tmo+1}\left[  T^{\tilde{m}_0+1}_{\tilde{m}_0}
> k-1
\right] & \leq &   \frac{3 (\gamma  \log_2 n)^{9/2}
}{\sigma^{10}(\gamma )^{1/2}(\log
n)^{\gamma-(\frac{12}{\kappa}+\frac{21}{2})}} .
\end{eqnarray}
Finally, using \ref{superdeltarp} and therefore choosing
$q=[q_n]$, where $q_n$ is given at the end of Definition
\ref{super},  we get that for all $\gamma
> \frac{12}{\kappa}+\frac{21}{2}$, $n>n_1$ and  $\alpha \in G_n$
\begin{eqnarray}
 \sum_{ q=[q_n] \leq k
\leq n-2} \alpha_{\tmo} \p^{\alpha}_{\tmo+1}\left[
T^{\tilde{m}_0+1}_{\tilde{m}_0}
> k
\right]  & \leq & \sum_{q=[n^{\delta_{r,r}}]+1 \leq k \leq n-2}
\alpha_{\tmo} \p^{\alpha}_{\tmo+1}\left[
T^{\tilde{m}_0+1}_{\tilde{m}_0}
> k
\right]  \\
& \leq & \frac{3 (  \log_2 n)^{9/2} }{\sigma^{10}(\gamma
)^{1/2}(\log n)^{\gamma-(\frac{12}{\kappa}+\frac{21}{2})}}
\label{4eq79}
\end{eqnarray}
Making similar computation for the sum on the right hand side of
\ref{eqtamere2} one get the same upper bound with $q=[q_n]$. Using
these estimates, \ref{eqtamere2}, \ref{eqtamere1},
\ref{thshypsecmom1} and the fact $\left\{\A_q^c, T^0_{\tmo} \leq
n\right\} \subset \bar{\A}_q^c $ we get the lemma taking $q=[q_n]$
 and $n_1''=n_1$.
\end{Pre}

\noindent We get Proposition \ref{ProS1} collecting the results of
Lemmata \ref{lempro1},  \ref{lempro2}, using \ref{6eq8} and taking
$n_0'=n_1' \vee n_1''$ and  $q=[q_n$].
\end{Pre}
\noindent \\

\begin{Pro}   \label{ProS2} There exists $h>0$, such that if \ref{hyp1} and  \ref{hyp0} hold
and for all $\kappa \in ]0,\kappa^+[$ \ref{hyp4} holds, for all
 $\gamma>0$ there exists $n_0\equiv
n_0(\gamma,\kappa)$ such that for all $n>n_0$ there exists $G_n
\subset \Omega_1$ with  $Q\left[G_n\right]\geq 1-h \left(( \log_3
n)( \log_2 n)^{-1}\right)^{1/2} $ and for all $\alpha \in G_n$

\begin{eqnarray}
\p^{\alpha}_0\left[\left|\frac{X_n}{(\log
n)^2}-m_{0}\right|>\delta_{q_n},\ \A_{q_n} \right]\leq
\frac{1}{(\log n)^{\gamma}} , \label{4.62}
\end{eqnarray}
$\delta_{q_n}= L_n (\log n)^{-2}$, $q_n$ and $L_n$ are given at
the of definition \ref{super}.
\end{Pro}
\begin{Pre}
Let us introduce the following stopping time $T_{\tmo}(q)=\inf
\left\{l \geq n-q,\ X_l= \tmo  \right\}$. We remark that $\A_q
\Leftrightarrow n-q \leq T_{\tmo}(q) \leq n$. Taking $q=[q_n]$, by
the strong Markov property we have
\begin{eqnarray}
 \p^{\alpha}_0\left[\left|\frac{X_n}{(\log
n)^2}-m_{0}\right|>\delta_{q_n},\ \A_{[q_n]}\right] &=&
\sum_{l=n-{[q_n]}}^n \p^{\alpha}_{\tmo}\left[
\left|\frac{X_{n-l}}{(\log n)^2}-m_{0}\right|>\delta_{q_n}
\right]\p^{\alpha}_0\left[T_{\tmo}({q_n})=l \right] \label{4.71}.
\end{eqnarray}
Therefore we get
\begin{eqnarray}  \p^{\alpha}_0\left[\left|\frac{X_n}{(\log
n)^2}-m_{0}\right|>\delta_{q_n},\ \A_{[q_n]} \right]  & \leq &
\sum_{l=0}^{q_n} \p^{\alpha}_{\tmo}\left[
T^{\tilde{m}_0}_{\tilde{m}_0+L_n}\wedge
T^{\tilde{m}_0}_{\tilde{m}_0-L_n} <
{q_n}-l\right]\p^{\alpha}_0\left[T_{\tmo}({q_n})=l \right]
\\
& \leq & \p^{\alpha}_{\tmo}
\left[T^{\tilde{m}_0}_{\tilde{m}_0+L_n}\wedge
T^{\tilde{m}_0}_{\tilde{m}_0-L_n} < {q_n} \right],
\label{ilest2heure}
\end{eqnarray}
Using Lemma \ref{lemsauve} we get \ref{4.62}.
\end{Pre}
\\
\noindent
Now we end the proof of theorem \ref{thSinai1}  \\
Assume \ref{hyp1},  \ref{hyp0} hold, let $\kappa \in ]0,\kappa^+[$
such that \ref{hyp4} hold, let us denote $\gamma_0 =
\frac{12}{\kappa}+\frac{21}{2}$, let $\gamma>\gamma_0$. Taking
$q=[q_n]$ and $\delta_q=L_n (\log n)^{-2}$ in \ref{ladereqs} we
obtain from Propositions \ref{ProS1} and \ref{ProS2} that there
exists $n_1 \equiv n_1(\kappa,\gamma)$ such that for all $n>n_1$
and  all $\alpha \in G_n $
\begin{eqnarray}
\qquad \p^{\alpha}_0\left[\left|\frac{X_n}{(\log
n)^2}-m_{0}\right|>\delta_{q_n}\right] \leq \frac{3( \log_2
n)^{9/2}}{\sigma^{10}(\gamma)^{1/2} (\log
n)^{\gamma-\gamma_0}}+\Er\left(\frac{1}{(\log
n)^{\gamma-(6/\kappa+4)}} \right), \label{eqProS11}
\end{eqnarray}
 Moreover we remark
 that one can find $n_2
> n_1$ such that for all $n>n_2$ we have $
\delta_{q_n} \equiv L_n (\log n)^{-2} \leq \gamma (1600)^2
  ( \log_2 n)^{9/2} (\log n)^{-1/2} $.
$\blacksquare$
\\

\appendix

 \begin{center} {\Large APPENDIX}
\end{center}

\section{ {\large Proof of the good properties for the environment
(Proposition \ref{profonda})} \label{sec6}}

In all this section we will use standard facts on sums of i.i.d.
random variables, these results are summarized in the Section B of
this appendix.

\subsubsection{Elementary results on the basic valley $\{\tM_0',\tmo,\tM_0\}$ }
\noindent We introduce the following stopping times, for $a>0$,
\begin{eqnarray}
& & U^+_a\equiv U^+_a(S_j^n,j\in \N) = \left\{\begin{array}{l}  \inf  \{m \in \N^*,\ S_m^n \geq a \} , \label{def1}  \\
 + \infty \textrm{, if such a } m \textrm{ does not exist.} \end{array} \right. \label{4.3} \\
& & U^-_a\equiv U^-_a(S_j^n,j\in \N)= \left\{\begin{array}{l}  \inf  \{m \in \N^*,\ S_m^n \leq -a \} , \label{4.3bbc}  \\
 + \infty \textrm{, if such a } m \textrm{ does not exist.}
\end{array} \right.
\end{eqnarray}

\noindent \textbf{Proof of lemma \ref{moexiste}}  To prove this
lemma it is enough to prove that the valley
$\{U^-_{1+\gamma(n)},\tm,U^+_{1+\gamma(n)}\}$  satisfies the three
properties of Definition \ref{thdefval1b} with a probability very
near 1. Let $\kappa \in ]0,\kappa^+[$, and $\gamma>0$.
 \noindent By definition of $U^-_{1+\gamma(n)}$ and
 $U^+_{1+\gamma(n)}$,
 $\{U^-_{1+\gamma(n)},\tm,U^+_{1+\gamma(n)}\}$ satisfies the two first properties of Definition \ref{thdefval1b}. We are left with the third property.
Assume $\tm>0$, we remark that $ S^{n}_{U^-_{1+\gamma(n)}}-\max_{0
\leq t \leq m}\left(S^{n}_t\right)\leq \gamma(n) \Rightarrow
\max_{0 \leq t \leq m}\left(S^{n}_t\right)\geq 1 $ moreover $
\max_{0 \leq t \leq \tm}\left(S^{n}_t\right) \leq 1+\gamma(n) \ .$
Therefore
\begin{eqnarray}
Q\left[S^{n}_{U^-_{1+\gamma(n)}}-\max_{0 \leq t \leq
m}\left(S^{n}_t\right)\leq \gamma(n) \right] & \leq  & Q\left[1
\leq \max_{0 \leq t \leq \tm}\left(S^{n}_t\right) \leq
1+\gamma(n)\right] \label{4127s} .
\end{eqnarray}
Using \ref{sym} and Lemma \ref{lem101}, it is easy to prove that
there exists $n_1\equiv n_1(\gamma,\sigma,
\E\left[|\epsilon_0|^3\right])$ such that for all $n>n_1$
\begin{eqnarray}
Q\left[S^{n}_{\tm}\leq -\gamma(n) \right]\geq 1-\frac{ \log_2
n}{\log n}\left(\gamma+\Er\left(\frac{ 1 }{\log_2 n}\right)\right)
\ . \label{exi1b}
\end{eqnarray}
Let us denote $\A=\{1 \leq \max_{0 \leq t \leq
\tm}\left(S^{n}_t\right) \leq 1+\gamma(n),\ S^{n}_{\tm}\leq
-\gamma(n)\}$, by \ref{4127s} and \ref{exi1b} we have
\begin{eqnarray}
Q\left[S^{n}_{U^-_{1+\gamma(n)}}-\max_{0 \leq t \leq
m}\left(S^{n}_t\right)\leq \gamma(n) \right] \leq   Q\left[\A
\right]+\frac{ \log_2 n}{\log n}\left(\gamma+\Er\left(\frac{ 1
}{\log_2 n}\right)\right)\ . \label{exi1}
\end{eqnarray}
Let us define
\begin{eqnarray}
W_{\gamma(n)}= \left\{\begin{array}{l}  \inf\{m \in \N^*,\ S_m^n \in [1,1+\gamma(n)]\} \ , \label{def1}  \\
 + \infty \textrm{, if such } m \textrm{ does not exist.} \end{array} \right.
\end{eqnarray}
Denote $\A'=\bigcup_{j>W_{\gamma(n)}}\left\{ S_j^n \leq
-\gamma(n),\ \bigcap_{k=W_{\gamma(n)}+1}^j\left\{ S^n_k <
1+\gamma(n)\right\} \right \}$, we have $\A \subset\A'$ so $
Q\left[\A \right]  \leq Q\left[\A'\right]$. Making a partition on
the values of $W_{\gamma(n)}$, using that $\{W_{\gamma(n)}=r \}
\Rightarrow \{S_r^n \in [1,1+\gamma(n)]\}$ and the strong Markov
property we get
\begin{eqnarray}
 Q\left[\A'\right] & \leq & \sup_{1-\gamma(n)  \leq x \leq 1
}\left( Q\left[U_{\gamma(n)+x}^-<
U_{1+\gamma(n)-x}^+\right]\right) \sum_{r=0}^{+\infty}
\int^{1+\gamma(n)}_{1}
Q\left[W_{\gamma(n)}=r,\ S_r^n \in dx \right] \\
& \leq & Q\left[U_{1}^-< U_{ 2 \gamma(n)}^+\right] \ .\label{exi3}
\end{eqnarray}
 Using lemma
\ref{lem101},  we get that there exists $n_2\equiv
n_2(\sigma,\E[|\epsilon_0|^3])$ such that for all $n>n_2$
\begin{eqnarray}
Q\left[U^-_{1}<U^+_{2\gamma(n)}\right] \leq \frac{2 \log_2 n}{\log
n}\left(\gamma+O\left(\frac{ 1 }{\log_2 n}\right)\right) \ .
\label{exi4}
\end{eqnarray}
Collecting what we did above and taking $n_0=n_1 \vee n_2$ we get
the lemma. $\blacksquare$

\noindent \\ \textbf{Proof of proposition \ref{8eq18}},

\noindent Let us prove \ref{2eq126}, noticing that $\tM_0 \leq
U_{1+\gamma(n)}^+$, and using remark \ref{sym}, for all $G>0$ we
get
\begin{eqnarray}
Q\left[\tM_0 >  ( \sigma^{-1} \log n)^2 \log_2 n \right] & \leq &
Q\left[U^+_{1+\gamma(n)} \wedge U_G^-
>  ( \sigma^{-1} \log n)^2 \right]+ Q\left[U_{1}^+ \geq U_G^-  \right] \
.\label{6.82}
\end{eqnarray}
Taking $G=\left(\frac{2 \log_2 n}{ h_1^2 \log_3  n}\right)^{1/2} $
with $h_1>0$ and using \ref{lem101eq1}, we get that there exists
$n_1 \equiv n_1(h_1,\sigma, \E_Q\left[|\epsilon_0|^3\right])$ such
that for all $n>n_1$
\begin{eqnarray}
Q\left[U_{1+\gamma(n)}^+ \wedge U_G^-
> E (\log n)^2 \right]  \leq  2q_1^{  \frac{h_1 }{16}\log_3  n } \ ,
\end{eqnarray}
where $q_1<0.7$. Choosing correctly the numerical constant $h_1$
we get for all $n>n_1$:
\begin{eqnarray}
Q\left[U_{1+\gamma(n)}^+ \wedge U_G^-
> ( \sigma^{-1} \log n)^2 \log_2 n \right] \leq  \frac{1}{ \log_2 n} \label{2.109} \ .
\end{eqnarray}
Taking $D = \log n $ in \ref{lem101eq2} we get for all $n>n_1 $
\begin{eqnarray}
 Q\left[U_{1+\gamma(n)}^+ \geq U_G^-  \right] \leq \frac{1}{G} +
 \Er \left(\frac{(\log_2 n)^{3/2}}{\log n}\right) .\label{2.110}
\end{eqnarray}
Using \ref{6.82}, \ref{2.109}, \ref{2.110} and the expression of
$G$ we get \ref{2eq126}, the proof of \ref{2eq126b} is similar.
$\blacksquare$

\noindent \\ We recall that for all $\kappa \in ]0,\kappa^+[$, $
C\equiv C(\kappa)= \E_Q\left[e^{\kappa \epsilon_0}\right] \vee
\E_Q\left[e^{-\kappa \epsilon_0}\right]< +\infty \label{C}$.

\noindent \\ \textbf{Proof of lemma \ref{interminib}.} \noindent
Denote
\begin{eqnarray}
A_0= \left\{ \tM_0 \geq ( \sigma^{-1} \log n)^2 \log_2 n,\tM_0'
\leq -( \sigma^{-1} \log n)^2 \log_2 n \right\} \label{A0}.
\end{eqnarray}
Let $u_n=\left[((\log n)^2)(65\sigma^2(\log_2 n))^{-1}\right]+1$
and $v_n$ a sequence such that $u_n\times v_n=[( \sigma^{-1}\log
n)^2 \log_2 n ]+1$. Using \ref{8eq18} we know that there exists
$n_0' \equiv n_0'\left(\epsilon,\sigma, \E_Q\left[
\left|\epsilon_0 \right|^3\right]\right) $ such that for all
$n>n_0'$
\begin{eqnarray}
Q\left[\tM_0-\tmo \leq u_n  \right]\leq Q\left[\tM_0-\tmo \leq
u_n,\ A_0\right]+h\left(\frac{ \log_3  n}{ \log_2 n}\right)^{1/2}
. \label{eqquipu4}
\end{eqnarray}
\noindent We recall that, in all this work, $h$ is a strictly
positive numerical constant that can grow from line to line if
needed. \noindent Let us denote $B_{n,\sigma}=\{-[( \sigma^{-1}
\log n)^2 \log_2 n]-1,[( \sigma^{-1} \log n)^2 \log_2 n],\cdots,[(
\sigma^{-1} \log n)^2 \log_2 n]+1\} $, by definition
$S_{\tM_0}-S_{\tm_0} \geq \log n$, so
\begin{eqnarray}
 Q\left[\tM_0-\tmo \leq u_n,\ A_0\right]
 & \leq & Q\left[\max_{m \in B_{n,\sigma} }\max_{m
\leq l \leq m+u_n }\max_{m \leq j \leq m+u_n
}\left(|S_l-S_j|\right) \geq \log n \right] . \label{5.49}
\end{eqnarray}
Making similar computations to the ones did in the proof of
 \ref{lem102eq2} we get that there exists $n_1
\equiv n_1 (\sigma, C, \kappa )$ such that for all $n>n_1$,
\begin{eqnarray}
Q\left[\max_{m \in B_{n,\sigma} }\max_{m \leq l \leq m+u_n
}\max_{m \leq j \leq m+u_n }\left(|S_l-S_j|\right) \geq \log n
\right] \leq \frac{4 \log_2 n }{ \sigma ^2 (\log
 n)^{1/33}} ,
\end{eqnarray}
using \ref{eqquipu4}, \ref{5.49}, \ref{eqquipu4} and taking
$n_0=n_0' \vee n_1$ we get \ref{intermini1b}. Similar computations
give \ref{intermini2b}. $\blacksquare$

\noindent \\ The following result is essential to the proof of the
other good properties.

\subsubsection{Minimal distance between the two points of one refinement (property \ref{superr})}

\begin{Lem} \label{lem103} There exists $h>0$ such that if \ref{hyp1bb}, \ref{hyp0} hold and for all $\kappa \in
]0,\kappa^+[$ \ref{hyp4} holds, for all $\gamma>0$ there exists
$n_0 \equiv n_0(\sigma
,\kappa,\E\left[|\epsilon_0|^3\right],C,\gamma)$ such that for all
 $n>n_0$
\begin{eqnarray}
Q\left[\bigcup_{i=1}^{r'}\left\{\tM_i'-\tm_i' \leq b_n \right\}
\right]\leq  h \left(\frac{\log_3 n}{\log_2 n} \right)^{1/2}+\Er
\left(\frac{ \log_2 n}{(\log n)^{1/33}}\right),
\label{eq2.144} \\
 Q\left[\bigcup_{i=1}^{r}\left\{\tM_i-\tm_i \leq b_n
\right\} \right]\leq h \left(\frac{\log_3 n}{\log_2 n}
\right)^{1/2}+\Er \left(\frac{ \log_2 n}{(\log n)^{1/33}}\right) .
\label{eq2.145}
\end{eqnarray}
$b_n$ is given in \ref{4.10b}, $\tM_{.}'$, $\tm_{.}'$ $\tM_{.}$
and $\tm_{.}$ have been defined Section \ref{OC}.
\end{Lem}

\begin{Rem}
This lemma shows that the distance between two points obtained by
 the operation of refinement is larger than $b_n$.
\end{Rem}

\begin{Pre}
Let $\kappa \in ]0,\kappa^+[$ and $\gamma>0$. Recalling
\ref{4.11b} and \ref{4.12b}, let us denote
\begin{eqnarray}
& & A_1= \bigcup_{i=1}^{r'}\left\{\tM_i'-\tm_i' \leq b_n \right\}
\\
& &
A_2=\bigcup_{l=-[k_n]-1}^{[k_n]+1}\bigcup_{j=l+[l_n]}^{[k_n]+1}
\left\{ \max_{(l+1)b_n \leq w < z \leq jb_n}\left(S_z-S_w \right)
\leq \max_{ lb_n \leq m \leq (j+1)b_n }\max_{m \leq u<v \leq m+b_n
} \left(S_v-S_u\right) \right\} .
\end{eqnarray}
Denoting $
C_1=\bigcap_{j=0}^{r'}\bigcup_{l=-[k_n]-1}^{[k_n]+1}\left\{\tM_j'
\in [lb_n,(l+1)b_n]\right\} $ and
$D_1=\bigcup_{i=1}^{r'}\bigcup_{l=-[k_n]-1}^{[k_n]+1}\Big\{\tM_i'-\tm_i'
\leq b_n,$ $\ \tM_{i-1}' \in [lb_n,(l+1)b_n]\Big\}$, it is clear
that $\{A_1,C_1\} \subset \{D_1\}$. Now denoting
$C_2=\bigcap_{i=0}^{r'-1}\left\{\tM_i' \leq \tmo-l_nb_n \right\}$
and
$D_2=\bigcup_{i=1}^{r'}\bigcup_{l=-[k_n]-1}^{[k_n]+1}\left\{\tM_i'-\tm_i'
\leq b_n,\ \tM_{i-1}' \in [lb_n,(l+1)b_n],\ \tM_{i-1}' \leq
\tmo-l_nb_n \right\}$, we easily get that $\{D_1,C_2\} \subset
D_2$. Finally denoting $
C_3=\bigcup_{l=-[k_n]-1}^{[k_n]+1}\left\{\tm_0 \in
[lb_n,(l+1)b_n]\right\} $, $D_3=
\bigcup_{i=1}^{r'}\bigcup_{l=-[k_n]-1}^{[k_n]+1}\bigcup_{j=l+[l_n]}^{[k_n]+1}$
$\left\{\tM_i'-\tm_i' \leq b_n,\ \right.$ $\left. \tM_{i-1}' \in
[lb_n,(l+1)b_n],\
    \tmo \in [b_nj,b_n(j+1)] \right\}$ and noticing that $\left\{ \tM_{i-1}'
\leq \tmo-l_nb_n, \ \tM_{i-1}' \in [lb_n,(l+1)b_n]\right\} \subset
\left\{\tmo \geq lb_n+l_nb_n \right\}$, we get that $\{D_2,C_3\}
\subset D_3$.
 Moreover if we
make a refinement of $\{\tM_{i-1},\tmo\}$, we get the points
$\tM_i'$ and  $\tm_i'$ such that $
S_{\tM_i'}-S_{\tm_i'}=\max_{\tM_{i-1}' \leq w < z \leq \tmo
}\left(S_z-S_w \right)$, so $D_3 \subset A_2$. Therefore we have :
\begin{eqnarray}
Q[A_1]\leq Q[A_2]+Q[C_1^c] +Q[C_2^{c}]+Q[C_3^c] \label{123}
\end{eqnarray}
It is easy to see that $\{C_1^c \subset A_0^{c}\}$, $\{C_1^c
\subset A_0^{c}\}$ and $C_2^c \subset \{\tmo-\tM_0' \geq (\log
n)^2(65 \sigma^2 \log_2 n)^{-1} \}$ so using Proposition
\ref{superint} and Lemma \ref{interminib} we have some upper
bounds for the three last probabilities of \ref{123}.

\noindent Now let us give an upper bound for $Q[A_2]$, first we
introduce the following event, let $s>0$
\begin{eqnarray}
\qquad A_3=\max_{-([k_n]+1)b_n \leq m \leq ([k_n]+1)b_n }\max_{m
\leq l \leq m+b_n }\max_{m \leq j \leq m+b_n
}\left(\left|S_l-S_j\right| \right)\leq g_n ,
\end{eqnarray}
where $g_n=((1+s) 32 \sigma^2 b_n \log k_n)^{1/2}$, we have
\begin{eqnarray}
Q\left[A_2\right]  \leq Q\left[A_2,A_3\right]+ Q\left[A_3^c\right]
\label{2eq223} .
\end{eqnarray}
Applying inequality \ref{lem102eq2}, (taking $[L]+1=([k_n]+1)b_n$
and $\log K= \log (k_n)$) we get that there exists $n_1 \equiv n_1
(\sigma , s, \kappa, \E\left[|\epsilon_0|^3,C\right]) $ such that
for all $n>n_1$
\begin{eqnarray}
 Q\left[A_3^c\right]\leq \frac{4b_n}{k_n^{\frac{s}{2}}} \label{es2}.
\end{eqnarray}
We are left to estimate $Q\left[A_2,A_3\right]$, we have
\begin{eqnarray}
 Q\left[A_2,A_3\right]
 & \leq &
 \sum_{i=-[k_n]-1}^{[k_n]+1}Q\left[\bigcup_{j=i+[l_n]}^{[k_n]+1}\left\{\max_{(i+1)
b_n\leq w < z \leq j b_n}\left(S_z-S_w\right) \leq  g_n
\right\}\right] .
\end{eqnarray}
We remark that the event $\left\{\max_{i b_n\leq w < z \leq j
b_n}\left(S_z-S_w\right) \leq  g_n \right\}$ is decreasing in $j$,
 so
\begin{eqnarray}
Q\left[A_2,A_3\right] &\leq&
\sum_{i=-[k_n]-1}^{[k_n]+1}Q\left[\max_{(i+1) b_n\leq w < z \leq
(i+[l_n])b_n}\left(S_z-S_w\right) \leq  g_n \right] .
\label{2en1eq2}
\end{eqnarray}
Denoting $(a_n, n \in \N^*)$ and $(d_n, n \in \N^*)$ two strictly
positive increasing sequence such that $[l_n]=d_n \times a_n$ we
get by independence
\begin{eqnarray}
Q\left[A_2,A_3\right] =2([k_n]+1)\left(Q\left[S_{a_nb_n} \leq g_n
 \right]\right)^{[d_n]-1} . \label{2en1eq3}
\end{eqnarray}
Now applying the Berry-Essen theorem to $Q\left[S_{a_nb_n} \leq
g_n \right]$ and choosing $d_n=- 2 \frac{(\log(k_n+2))}{(\log
(\int_1^{+\infty} e^{-x^2}/(2\pi)^{1/2}))}$, we obtain that there
exists $n_2 \equiv n_2(\sigma, \E_Q[|\epsilon_0|^3])$ such that
for all $n>n_2$
\begin{eqnarray}
 Q\left[A_2,A_3\right] &\leq& \frac{2}{k_n} . \label{ineq1}
\end{eqnarray}
Finally, taking $s=4$ and using \ref{2eq223}, \ref{es2} and
\ref{ineq1} we get that there exists $n_3 \equiv n_3(
\sigma,\kappa,\E_Q\left[|\epsilon_0|^{3}\right],C,\gamma) \geq n_1
\vee n_2 $ such that for all $n>n_3$
\begin{eqnarray}
Q\left[A_2\right] = \Er \left(\frac{\log_2 n }{\log n
}\right)^{1/2}  \label{es1}
\end{eqnarray}
Collecting \ref{123} and \ref{es1} we get \ref{eq2.144}. \noindent
Similar computations give \ref{eq2.145}. \\
\end{Pre}

\begin{Cor} \label{superrb}  There exists $h>0$ such that if \ref{hyp1bb}, \ref{hyp0} hold and for all $\kappa \in
]0,\kappa^+[$ \ref{hyp4} holds, for all $\gamma>0$ there exists
$n_0 \equiv n_0(\sigma ,\E\left[|\epsilon_0|^3\right],C,\gamma)$
such that for all $n>n_0$
\begin{eqnarray}
Q\left[r' \leq 2k_n+1 \right]\geq 1-h \left(\frac{ \log_3  n}{
\log_2 n}\right)^{1/2} -\Er \left(\frac{ \log_2 n}{(\log
n)^{1/33}}\right) \label{superr1b},
\\
Q\left[r \leq 2k_n+1 \right]\geq 1-h \left(\frac{ \log_3  n}{
\log_2 n}\right)^{1/2} -\Er\left(\frac{ \log_2 n }{(\log
n)^{1/33}}\right) \label{superr2b}.
\end{eqnarray}
$r$ and $r'$ have been defined section \ref{OC} and $k_n$ is given
 in \ref{4.11b}.
\end{Cor}

\begin{Pre}
This corollary is an easy consequence of lemma \ref{lem103}, the
proof is omitted.
\end{Pre}

\subsubsection{Minimal distance between two maximums (properties \ref{supereta} and \ref{superetap})}

\begin{Pro} \label{profondaeta} There exists $h>0$ such that if \ref{hyp1bb}, \ref{hyp0} hold and
for all $\kappa \in ]0,\kappa^+[$ \ref{hyp4} holds, there exists
$n_0 \equiv n_0(\sigma, \kappa
,\E\left[|\epsilon_0|^3\right],\E\left[\epsilon_0^4\right],C,\gamma)$
such that for all $n>n_0$
\begin{eqnarray}
& & Q\left[ \bigcap_{i=0}^{r-1} \left\{ \eta_{i,i+1} \geq
\gamma(n) \right\} \right]\geq1-h \left(\frac{ \log_3  n}{ \log_2
n}\right)^{1/2} - \Er\left(\frac{1}{\log_2 n}\right) , \label{prosupereta} \\
& & Q\left[ \bigcap_{i=0}^{r'-1} \left\{ \eta_{i,i+1}' \geq
\gamma(n) \right\} \right]\geq 1-h \left(\frac{ \log_3  n}{ \log_2
n}\right)^{1/2} - \Er\left(\frac{1}{\log_2 n}\right) \
,\label{prosuperetap}
\end{eqnarray}
where $\gamma(n)$ is given a the end of Definition \ref{super},
$\eta_{.,.}$ and $\eta_{.,.}'$ are given in \ref{notsuper}.
\end{Pro}

\begin{Pre}
  \textbf{Let us prove \ref{prosuperetap}} \\
To prove this proposition we will use the lemma \ref{lem103}. Let
$n>3$, and $\gamma>0$, we recall the following notations $
b_n=\left[(\gamma )^{1/2} (\log n \log_2 n )^{3/2} \right]+1$, $
k_n=((\sigma^{-1}\log n)^2 \log_2 n)/b_n$. Let us denote
\begin{eqnarray}
&& A=\bigcap_{i=0}^{r'}\left\{-(\sigma^{-1}\log n)^2 \log_2 n \leq
\tM_i' \leq (\sigma^{-1}\log n)^2 \log_2 n \right\} ,
\\
&& A_1=\bigcup_{i=1}^{r'}\bigcup_{j=-[k_n]-1}^{[k_n]+1}\left\{m_i'
\in
[b_nj,b_n(j+1)],\ M_i' \in [b_nj,b_n(j+1)]\right\} ,\\
&& A_2=\bigcup_{i=1}^{r'}\bigcup_{j=-[k_n]-1}^{[k_n]+1}\left\{M_i'
\in [b_nj,b_n(j+1)],\ M_{i+1}'\in [b_nj,b_n(j+1)]\right\} , \\
&& A_3=\bigcup_{i=0}^{r'-1}\left\{0 \leq \eta_{i,i+1}' \leq
\gamma(n) \right\} .
\end{eqnarray}
We have $Q[A_3] \leq Q[A_3,A_1^c,A]+Q[A_1]+Q[A^c]$, moreover $A
\subset A_0$ (see \ref{A0}) and $ A_1 \subset
\bigcup_{i=1}^{r'}\left\{\tM_i'-\tm_i' \leq b_n \right\} $,
therefore using Lemma \ref{superint} and the inequality
\ref{eq2.144} we get that there exists $h>0$ and $n_1$ such that
for all $n>n_1$, $Q[A_3] \leq Q[A_3,A_1^c,A]+h ((\log_2 n)/(\log
n))^{1/2} $. Let us denote $L_{i,j}(n)=\max_{b_ni \leq k \leq
b_n(i+1)}\left(S^n_k\right)-\max_{b_nj \leq l \leq
b_n(j+1)}\left(S^n_l\right)$, define
\begin{eqnarray}
A_4=\bigcup_{i=-[k_n]-1}^{[k_n]+1}\bigcup_{j=i+1}^{[k_n]+1}\left\{0
\leq L_{i,j}(n) \leq \gamma(n)\right\} ,
\end{eqnarray}
by definition of the refinements we have  $\tM_i'<\tM_{i+1}'
\textrm{ and  } S_{\tM_i'}>S_{\tM_{i+1}'}, \forall
 i \  0\leq i \leq r'-1$,
therefore $\{A_3,\ A_2^c,\ A\} \subset A_4$ then $Q\left[A_3,\
A_2^c,\ A\right]  \leq  Q\left[A_4 \right] $. Finally, we get that
for all $n>n_3$
\begin{eqnarray}
Q[A_3] \leq Q[A_4]+h ((\log_2 n)/(\log n))^{1/2} \label{Leqeq}
\end{eqnarray}
 Denoting
\begin{eqnarray}
& &
A_5=\bigcup_{i=-[k_n]-1}^{[k_n]+1}\bigcup_{j=i+2}^{[k_n]+1}\left\{0\leq
L_{i,j}(n) \leq \gamma(n) \right\},  \label{superetaeq2}\\
& & A_6=\bigcup_{i=-[k_n]-1}^{[k_n]+1}\left\{0 \leq L_{i,i+1}(n)
\leq \gamma(n) \right\}
 . \label{superetaeq3}
\end{eqnarray}
we have that
\begin{eqnarray}
Q\left[A_4 \right] = Q\left[A_5 \right]+Q\left[A_6 \right].
\label{Leqeqeq}
\end{eqnarray}
Now we estimate the two probability $Q\left[A_5 \right]$ and
$Q\left[A_6 \right]$ in (respectively) lemma \ref{superetalem1}
and \ref{superetalem2}. For the proof of these lemmata we have
used the paper in preparation of \cite{Picco}.
\begin{Lem} \label{superetalem1} Assume \ref{hyp1}, \ref{hyp0} and  \ref{hyp4}, for all $\gamma>0$
there exists $n_0' \equiv
n_0'\left(\sigma,\gamma,\E\left[\epsilon_0 ^4 \right] \right)$
such that for all $n>n_0'$

\begin{eqnarray}
Q\left[A_5 \right] & \leq & 10
\left(\frac{\pi}{\sigma^2}\right)^{1/2}\frac{\gamma \log_2
n}{(b_n)^{1/2}}\left([k_n]+1\right)^{3/2}
\end{eqnarray}
where $k_n$ is given by \ref{4.11b}, $b_n$ by  \ref{4.10b}.
\end{Lem}
\begin{Pre} We have
\begin{eqnarray}
Q\left[A_5 \right]  \leq
\sum_{i=-[k_n]-1}^{[k_n]+1}\sum_{j=i+2}^{[k_n]+1}Q\left[\left\{0
\leq L_{i,j}(n) \leq \gamma(n) \right\} \right] .
\label{superetaeq4}
\end{eqnarray}
Now we give an upper bound for
$\sum_{i=0}^{[k_n]+1}\sum_{j=i+2}^{[k_n]+1}Q\left[\left\{0 \leq
L_{i,j}(n) \leq \gamma(n) \right\} \right] $. Denoting
$Z_{i+1,j}(n)=-\sum_{l=b_n(i+1)+1}^{b_n j}\epsilon_l$ and
$Y=-\min_{ i b_n \leq k \leq (i+1)b_n
}\sum_{m=k}^{(i+1)b_n}\epsilon_m-\max_{ j b_n+1 \leq k \leq
(j+1)b_n }\sum_{m=jb_n+1}^{k}\epsilon_m$, it is easy to see that
for all $i \geq 0$, $L_{i,j}(n)  = (Z_{i+1,j}(n)+Y)/(\log n)$.
\noindent Therefore we have
\begin{eqnarray}
 Q\left[0 \leq L_{i,j}(n) \leq \gamma(n) \right]  & =&\int_{\R}
Q\left[ 0 \leq Z_{i+1,j}(n)-y \leq \gamma(n) \log n,\ Y \in dy
\right] . \label{3eq254b}
\end{eqnarray}
$Z_{i+1,j}(n)$ and Y are independent so
\begin{eqnarray}
 \int_{\R} Q\left[ 0 \leq Z_{i+1,j}(n)-y \leq \gamma(n) \log n,\ Y \in dy \right]  & \leq &
\sup_y\left(Q\left[ y \leq Z_{i+1,j}(n) \leq \gamma(n)\log
n+y\right]\right) . \label{derder}
\end{eqnarray}
To estimate this last term we use the following concentration
inequality
 (see \cite{Lecam} pages 401-413)
\begin{eqnarray}
\sup_y\left(Q\left[ y \leq Z_{i+1,j}(n) \leq \gamma(n)\log
n+y\right]\right) \leq \frac{2(\pi)^{1/2}}{Z} \label{3eq258b} ,
\end{eqnarray}
where $Z^2 \equiv Z^2(\gamma(n))=\sum_{l=1}^{b_n (j-i-1)}\E\left[1
\wedge H_s^2\right]$, $H_s=\frac{\epsilon^s_l}{\gamma(n)\log n}$
and $\epsilon^s_l=\epsilon_l-\epsilon_s'$, $\epsilon_l'$
 is independent and identically distributed to $\epsilon_l$. We have
$ \E\left[1 \wedge (H_s)^2 \right] \geq (\gamma(n)\log
n)^{-2}\E\left[\left(\epsilon^s_l\right)^2\un_{1>H_s}\right] $.
Noticing that
$\E\left[\left(\epsilon^s_l\right)^2\un_{1>H_s}\right] =
\E\left[\left(\epsilon^s_l\right)^2\right]-\E\left[\left(\epsilon^s_l\right)^2\un_{1\leq
H_s}\right]$ we get by Schwarz inequality and Markov inequality
\begin{eqnarray}
\E\left[\left(\epsilon^s_l\right)^2\un_{1>H_s} \right] &\geq  & 2
\sigma^2 - \left(\E\left[\left(\epsilon^s_l\right)^4
\right]^{1/2}(2 \sigma^2)^{1/2} \right)(\gamma \log_2 n)^{-1} .
\end{eqnarray}
 We deduce that there exists $n_0'
\equiv n_0'\left(\sigma,\gamma,\E\left[\epsilon_0 ^4 \right]
\right)$ such that for all $n>n_0'$, $\E\left[1 \wedge
\left(H_s\right)^2\right] \geq 3\sigma^2 /(2(\gamma(n) \log n)^
2)$, therefore for all $n>n_0'$
\begin{eqnarray}
Z \geq \sqrt{\frac{3}{2}\sigma^2} \frac{\sqrt{b_n(j-i-1)}}{
\gamma(n)\log n} \label{Leq}.
\end{eqnarray}
Inserting \ref{Leq} in \ref{3eq258b} and using \ref{derder} and
\ref{3eq254b} we obtain for all $n
> n_0'$
\begin{eqnarray}
 Q\left[0 \leq L_{i,j}(n) \leq \gamma(n) \right]
& \leq &
\left(\frac{8\pi}{3\sigma^2}\right)^{1/2}\frac{\gamma(n)\log
n}{(b_n)^{1/2}\left(j-i-1\right)^{1/2}} \label{3eq265} .
\end{eqnarray}
Therefore, using \ref{3eq265} for all $n>n_0'$ we have
\begin{eqnarray}
 \sum_{i=0}^{[k_n]+1}\sum_{j=i+2}^{[k_n]+1}Q\left[\left\{0 \leq
L_{i,j}(n) \leq \gamma(n) \right\} \right]  &\leq &
\frac{5}{2}\left(\frac{\pi}{\sigma^2}\right)^{1/2}\frac{\gamma
\log_2 n}{(b_n)^{1/2}}\left([k_n]+1\right)^{3/2} \label{3eq267}
\end{eqnarray}
Making similar computations for the case $i < 0$ we get a similar
result, so we get lemma \ref{superetalem1}.
\end{Pre}
\\

\noindent \textbf{Constraint on $k_n$ and  $b_n$ }
\label{supercontrainte} Now we can justify the choice for $b_n$
and $k_n$, recalling that $k_n\times b_n=(\sigma^{-1}\log n)^2
\log_2 n$ we want that
\begin{eqnarray}
 \left(\frac{\pi}{\sigma^2}\right)^{1/2}\frac{\gamma \log_2
n}{(b_n)^{1/2}}\left([k_n]+1\right)^{3/2} \label{superetaeq9} ,
\end{eqnarray}
be close to $0$ but  $b_n$ small.
Using that $b_n=\left[(\gamma )^{1/2} (\log n \log_2 n
)^{3/2}\right]+1$, we get that there exists $h_1\equiv
h_1(\sigma,\gamma)
>0$ and $n_2$ such that for all $n>n_2$,
\begin{eqnarray}
10 \left(\frac{\pi}{\sigma^2}\right)^{1/2}\frac{\gamma \log_2
n}{(b_n)^{1/2}}\left([k_n]+1\right)^{3/2} & \leq & h_1
\left(\frac{1}{ \log_2 n}\right)^{1/2}. \label{3eq272cc}
\end{eqnarray}
 So using
\ref{3eq272cc} and lemma \ref{superetalem1}, we get that there
exists $n_1' \equiv n_1'(\sigma,\gamma,\E[\epsilon_0^4]) \geq n_0'
\vee n_2 $ such that for all $n>n_1'$
\begin{eqnarray}
& &
Q\left[\bigcup_{i=-[k_n]-1}^{[k_n]+1}\bigcup_{j=i+2}^{[k_n]+1}\left\{\max_{b_ni
\leq k \leq b_n(i+1)}\left(S^n_k\right)-\max_{b_nj \leq l \leq
b_n(j+1)}\left(S^n_l\right)\leq \gamma(n) \right\} \right] \leq
h_1 \left(\frac{1}{ \log_2 n}\right)^{1/2} .\label{superetaeq9bis}
\end{eqnarray}
Now we prove the following lemma
\begin{Lem} \label{superetalem2} Assume \ref{hyp1}, \ref{hyp0} hold and for all $\kappa \in ]0,\kappa^+[$
  \ref{hyp4} holds, for all $\gamma>0$
 there exists $n_0'' \equiv n_0''(\sigma ,\E\left[|\epsilon_0|^3\right],$ $\E\left[\epsilon_0^4\right],C,\gamma)$
 such that for all $n>n_0''$
\begin{eqnarray}
Q\left[A_6 \right] & \leq & \frac{(2[k_n]+3)(\log_2
n)^{5/2}}{(b_n)^{1/2}} \left(2\gamma+ \left(\frac{16 \pi
}{3\sigma^2}\right)\frac{\gamma }{\sigma (\log_2 n)^{3/2}}\right).
\end{eqnarray}
\end{Lem}

\begin{Pre} We have
\begin{eqnarray}
Q\left[\bigcup_{i=-[k_n]-1}^{[k_n]+1}\left\{0\leq L_{i,i+1}(n)
\leq \gamma(n)\right\}  \right]  &\leq&
\sum_{i=-[k_n]-1}^{[k_n]+1} Q\left[0 \leq L_{i,i+1}(n)\leq
\gamma(n) \right] \label{2eq328} .
\end{eqnarray}
Using the fact that we can write $\max_{b_n(i+1) \leq l \leq
b_n(i+2)}\left(S^n_l\right)=X+\max_{b_n(i+1)+1 \leq l \leq
b_n(i+2)}\left( \sum_{l=b_n(i+1)}^{l}\right)$ with $ X \in
\sigma\left( \epsilon_1, \cdots, \epsilon_{b_n(i+1)} \right)$ and
$Y\equiv \max_{b_ni \leq k \leq b_n(i+1)}\left(S^n_k\right) \in
\sigma\left( \epsilon_1, \cdots, \epsilon_{b_n(i+1)} \right)$ we
easily get by independence that
\begin{eqnarray}
 Q\left[0 \leq L_{i,i+1}(n) \leq \gamma(n)  \right]  & \leq &
\sup_{x}\left(Q\left[x \leq \max_{1 \leq k \leq b_n
}\left(S_k^n\right)\leq x+\gamma(n)  \right]\right) ,
\end{eqnarray}
 replacing this in \ref{2eq328}, we get
\begin{eqnarray}
& & Q\left[\bigcup_{i=-[k_n]-1}^{[k_n]+1}\left\{0\leq L_{i,i+1}(n)
\leq \gamma(n) \right\} \right]   \leq  (2[k_n]+3)
\sup_{x}\left(Q\left[x \leq \max_{1 \leq k \leq b_n
}\left(S_k^n\right)\leq x+\gamma(n)  \right]\right) .
\label{3eq274b}
\end{eqnarray}
To estimate $ \sup_{x}\left(Q\left[x \leq \max_{1 \leq k \leq b_n
}\left(S_k^n\right)\leq x+\gamma(n) \right]\right)$ we remark that
\begin{eqnarray}
Q\left[x \leq \max_{1 \leq k \leq b_n }\left(S^n_k\right)\leq
x+\gamma(n) \right] &=& Q\left[U^{+}_x \leq b_n \leq
U^{+}_{x+\gamma(n)}
\right] \\
&=& Q\left[U^{+}_x \leq \frac{b_n}{2},\ U^{+}_{x+\gamma(n)}\geq
b_n \right] \label{superetaeq10}
\\
&+& Q\left[\frac{b_n}{2} < U^{+}_x \leq b_n \leq
U^{+}_{x+\gamma(n)} \right] . \label{superetaeq11}
\end{eqnarray}
We have to estimate the two probability in \ref{superetaeq10} and
\ref{superetaeq11}. We begin with \ref{superetaeq11}, we remark
that
\begin{eqnarray}
\frac{b_n}{2} < U^{+}_x \leq b_n \leq U^{+}_{x+\gamma(n)}
\Rightarrow x \leq \max_{b_n/2 \leq k \leq b_n
}\left(S_k^n\right)\leq x+\gamma(n) ,
\end{eqnarray}
from this we deduce by the concentration inequality  (see
equations \ref{3eq258b} to \ref{3eq265}) that there exists $n_3
\equiv n_3(\sigma,\E\left[\epsilon_0^4\right])$ such that for all
$n>n_3$
\begin{eqnarray}
 Q\left[\frac{b_n}{2} < U^{+}_x
\leq b_n \leq U^{+}_{x+\gamma(n)} \right] & \leq &
\sup_y\left(Q\left[y \leq S_{b_n/2}^n  \leq y+\gamma(n)
\right]\right)  \leq  \left(\frac{ 16 \pi
}{3\sigma^2}\right)^{1/2}\frac{\gamma \log_2 n}{(b_n)^{1/2}} .
\label{superetaeq12}
\end{eqnarray}
 Now we estimate the probability in \ref{superetaeq10}, by the strong  Markov property we have
\begin{eqnarray}
 Q\left[U^+_x \leq \frac{b_n}{2},\ U^+_{x+\gamma(n)}\geq b_n \right]
&=& \sum_{l=0}^{b_n/2}\int_{x}^{x+\gamma(n)} Q\left[U^+_x
=l,S_l\in dy \right]Q\left[\ U^+_{x+\gamma(n)-y}\geq b_n-l \right]
,
\end{eqnarray}
moreover $x-y \leq 0$, therefore $Q\left[\ U^+_{x+\gamma(n)-y}\geq
b_n-l \right] \leq Q\left[\ U^+_{\gamma(n)}\geq b_n-l \right]$, so
we get
\begin{eqnarray}
& & Q\left[U^+_x \leq \frac{b_n}{2},\ U^+_{x+\gamma(n)}\geq b_n
\right] \leq Q\left[\ U^+_{\gamma(n)}\geq b_n/2 \right] .
\end{eqnarray}
To estimate this probability we use remark \ref{symb} and  lemma
\ref{lem101} (taking $c =\frac{\gamma \log_2 n}{\log n} $, $a
=\frac{(b_n)^{1/2}}{\log n (\log_2 n)^{3/2}}$, $L=b_n/2$ and
$D=\log n$), we get that there exists $n_4$ such that for all
$n>n_4$
\begin{eqnarray}
Q\left[\ U^+_{\gamma(n)}\geq b_n/2 \right] &\leq&
\frac{2\gamma(\log_2 n)^{5/2}}{(b_n)^{1/2}} . \label{superetaeq13}
\end{eqnarray}
 Inserting
 \ref{superetaeq12} and   \ref{superetaeq13} in (respectively) \ref{superetaeq10}
 and  \ref{superetaeq11} and using
 \ref{3eq274b} we get for all $n>n_4$
\begin{eqnarray}
&&  Q\left[\bigcup_{i=-[k_n]-1}^{[k_n]+1}\left\{0\leq L_{i,i+1}(n)
\leq \gamma(n) \right\} \right]  \leq \frac{(2[k_n]+3)(\log_2
n)^{5/2}}{(b_n)^{1/2}} \left(2\gamma+ \left(\frac{16 \pi
}{3\sigma^2}\right)\frac{\gamma }{\sigma (\log_2 n)^{3/2}}\right)
,
\end{eqnarray}
taking  $n_0''=n_3 \vee n_4$ we get Lemma \ref{superetalem2} .
\end{Pre}

\noindent
 Recalling \ref{4.10b} and \ref{4.11b}  we get from Lemma \ref{superetalem2}, that
for all $\kappa \in]0,\kappa^+[$, $\gamma>0$
 there exists $n_1'' \equiv n_1''(\sigma,\kappa ,\E\left[|\epsilon_0|^3\right],$ $\E\left[\epsilon_0^4\right],C,\gamma) \geq n_0'' $
 such that for all $n>n_1''$
\begin{eqnarray}
Q\left[\bigcup_{i=-[k_n]-1}^{[k_n]+1}\left\{0\leq L_{i,i+1}(n)
\leq \gamma(n) \right\}  \right]  &=&  \Er\left(\frac{(\log_2
n)^{1+3/4}}{ (\log n)^{1/4}}\right) . \label{3eq292b}
\end{eqnarray}
To end the proof of Proposition \ref{profondaeta}, we collect
\ref{3eq292b}, \ref{superetaeq9bis}, \ref{Leqeqeq}, and finally
\ref{Leqeq}, and we take $n_0= n_1 \vee n_1'  \vee n_1''$. We get
\ref{prosupereta} with similar computations.
\end{Pre}

\subsubsection{Distance minimal between the maximum and the minimum of one refinement (properties \ref{superdelta} and \ref{superdeltap})}

\begin{Pro} \label{profondadelta} There exists $h>0$ such that if \ref{hyp1}, \ref{hyp0} hold and for all $\kappa \in]0,\kappa^+[$ \ref{hyp4}
holds, for all $\gamma>0$
 there exists $n_0 \equiv n_0(\sigma ,\E\left[|\epsilon_0|^3\right],\E\left[\epsilon_0^4\right],C,\gamma)$ such that for all $n>n_0$
\begin{eqnarray}
& & Q\left[ \bigcap_{i=0}^{r-1} \left\{ \delta_{i+1,i+1} \geq
\gamma(n) \right\} \right]\geq 1- h \left(\frac{ \log_3  n}{
\log_2 n}\right)^{1/2} - \Er\left(\frac{1}{\log_2 n}\right) ,
 \label{prosuperdelta} \\
&& Q\left[ \bigcap_{i=0}^{r'-1} \left\{ \delta_{i+1,i+1}' \geq
\gamma(n) \right\} \right] \geq 1-h \left(\frac{ \log_3  n}{
\log_2 n}\right)^{1/2} - \Er\left(\frac{1}{\log_2 n}\right) ,
 \label{prosuperdeltap}
\end{eqnarray}
 where $\gamma(n)$ is given at the end of Definition \ref{super}, $\delta_{.,.}$ and $\delta_{.,.}'$ are given in \ref{notsuper}.
\end{Pro}

\begin{Pre}
First we remark that by construction the event $\{\delta_{i+1,i+1}
\geq \gamma(n)\}$ decrease in $i$, so $
 Q\left[ \bigcap_{i=0}^{r-1} \left\{ \delta_{i+1,i+1} \geq
\gamma(n) \right\} \right]  =  Q\left[  \delta_{r,r} \geq
\gamma(n)  \right]$, then we use the same method used to prove
Proposition \ref{profondaeta}.
\end{Pre}
\\

\subsubsection{Minimal distance between a minimum and $S_{\tmo}$ (properties \ref{supermu} and \ref{supermup})}

\begin{Pro} \label{profondamu} There exists $h>0$ such that if \ref{hyp1}, \ref{hyp0} hold and for all $\kappa \in]0,\kappa^+[$ \ref{hyp4}
holds, for all $\gamma>0$
 there exists $n_0 \equiv n_0(\sigma ,\E\left[|\epsilon_0|^3\right],\E\left[\epsilon_0^4\right],C,\gamma)$ such that for all $n>n_0$
\begin{eqnarray}
& & Q\left[ \bigcap_{i=0}^{r-1} \left\{ \mu_{i+1,0} \geq \gamma(n)
\right\} \right]\geq 1-h \left(\frac{ \log_3  n}{ \log_2
n}\right)^{1/2} -  \Er\left(\frac{1}{\log_2 n}\right) ,
 \label{prosupermu1} \\
&& Q\left[ \bigcap_{i=0}^{r'-1} \left\{ \mu_{i+1,0}' \geq
\gamma(n) \right\} \right] \geq 1-h \left(\frac{ \log_3  n}{
\log_2 n}\right)^{1/2} -  \Er\left(\frac{1}{\log_2 n}\right) ,
 \label{prosupermu2}
\end{eqnarray}
where $\gamma(n)$ is given at the end of Definition \ref{super},
$\mu_{.,.}$ and $\mu_{.,.}'$ are given in \ref{notsuper}.
\end{Pro}

\noindent The proof of this proposition is similar to the proof of
Proposition \ref{profondaeta} and is omitted.

\subsubsection{Control of the first and the last refinement (properties \ref{superdelta1}, \ref{superdelta1p}, \ref{superdeltar} and \ref{superdeltarp})}

\begin{Pro} \label{profondaor} There exists $h>0$ such that if \ref{hyp1}, \ref{hyp0} hold and
for all $\kappa \in]0,\kappa^+[$ \ref{hyp4} holds , for all
 $\gamma>0$
 there exists $n_0 \equiv n_0(\sigma ,\E\left[|\epsilon_0|^3\right],\E\left[\epsilon_0^4\right],C,\gamma)$ such that for all $n>n_0$
\begin{eqnarray}
& &  Q\left[\delta_{1,1} \leq 1-\gamma(n)\right] \geq 1-h
\left(\frac{ \log_3  n}{ \log_2 n}\right)^{1/2}
- \Er\left(\frac{1}{\log_2 n}\right) \label{prosuperdelta1} , \\
& &  Q\left[\delta'_{1,1} \leq 1-\gamma(n)\right] \geq 1-h
\left(\frac{ \log_3  n}{ \log_2 n}\right)^{1/2}
- \Er\left(\frac{1}{\log_2 n}\right)  \label{prosuperdelta1p} , \\
& & Q\left[\delta_{r,r}  \leq (\log(q_n))(\log n)^{-1}\right]\geq
1 -h \left(\frac{ \log_3  n}{ \log_2 n}\right)^{1/2}
-\Er\left(\frac{ (\log_2
n)^{11/2}}{(\log n)^{1/66}}\right) \label{prosuperdeltar} , \\
& & Q\left[ \delta_{r',r'}' \leq (\log(q_n))(\log n)^{-1}
\right]\geq 1 -h \left(\frac{ \log_3  n}{ \log_2 n}\right)^{1/2}
-\Er\left(\frac{ (\log_2 n)^{11/2}}{(\log n)^{1/66}}\right) ,
\end{eqnarray}
 where $\gamma(n)$ and $q_n$ are given at the end of Definition \ref{super}.
\end{Pro}

\begin{Pre}
 \noindent \textbf{Let us prove \ref{prosuperdelta1}}, by construction $\delta_{1,1} \leq 1+\gamma(n)
$. So we have to prove that the event $ -\gamma(n)  \leq
\delta_{1,1}-1  \leq \gamma(n)$ has a probability very near 0, to
do this make we make use similar computations used to prove
Proposition
\ref{profondaeta}. A similar remark work for \ref{prosuperdelta1p}. \\
\textbf{Let us prove \ref{prosuperdeltar}}, by construction we
have
\begin{eqnarray}
& & \tM_0' \leq \tM_r \leq \tM_0 \label{4eq305} , \\
 & & \tM_r-\tm_0 \leq l_n \times b_n \label{4eq306} .
\end{eqnarray}
Using \ref{4eq305} and proposition \ref{8eq18}, we know that
 there exists $n_1 \equiv n_1\left(\sigma, \E\left[
\left|\epsilon_0 \right|^3\right]\right) $ such that for all
$n>n_1$

\begin{eqnarray}
Q\left[ -(\sigma^{-1} \log n)^2 \log_2 n \leq \tM_r \leq
(\sigma^{-1} \log n)^2 \log_2 n\right] \geq 1- h \left((\log_3 n)
(\log_2 n)^{-1}\right)^{1/2}
 \label{4eq310}
\end{eqnarray}
 Let us make the following chopping $\left[(\sigma^{-1} \log
n)^2 \log_2 n+1\right]=b_n'\times k_n'$ with $b_n'=[l_n\times
b_n]+1$,  we have $\delta_{r,0} \geq \delta_{r,r}$, therefore,
denoting $L'(n)=\max_{-b_n'\times k_n' \leq m \leq b_n'\times k_n'
}\max_{m \leq j \leq m+b_n' }\max_{m \leq l \leq m+b_n'
}\left(\left|S_l^n-S_j^n\right| \right) $
\begin{eqnarray}
  \left\{ \begin{array}{l}
  -(\sigma^{-1} \log n)^2 \log_2 n \leq \tM_r \leq (\sigma^{-1} \log n)^2 \log_2 n  \\
\textrm{and  }   \tmo-\tM_r
 \leq l_n \times b_n .
\end{array} \right\}
 \Rightarrow \delta_{r,r} \leq \delta_{r,0} \leq L'(n).
\end{eqnarray}
From this and \ref{4eq310} we deduce that for all $n>n_1$ we have
\begin{eqnarray}
  Q\left[\delta_{r,r} \leq L'(n) \right]
 \geq 1- h \left((\log_3 n) (\log_2 n)^{-1}\right)^{1/2}
 . \label{4eq311}
\end{eqnarray}
Using \ref{lem102eq2} (with $K=k_n'$, $[L]+1=[(\sigma^{-1} \log
n)^2 \log_2 n]+1$, $B=b_n'$ and $s=4$) one can check that that
there exists $n_2 \equiv n_2\left(\sigma,s, \kappa,\E\left[
\left|\epsilon_0 \right|^3\right],C\right) $ such that for all
$n>n_2$
\begin{eqnarray}
  Q\left[L'(n) > ((1+s)32
\sigma^2 b_n' \log k_n')^{1/2} \right]
 &=&  \Er\left(\frac{ (\log_2
n)^{11/2}}{(\log n)^{1/66}}\right), \label{4eq313}
\end{eqnarray}
Using  \ref{4eq311} and  \ref{4eq313} we get that for all $n>n_2$
\begin{eqnarray}
Q\left[\delta_{r,r} (\log n) \leq (160 \sigma^2 b_n' \log
k_n')^{1/2}\right] \geq 1- h \left((\log_3 n) (\log_2
n)^{-1}\right)^{1/2}
 -\Er\left(\frac{ (\log_2
n)^{11/2}}{(\log n)^{1/66}}\right) .
\end{eqnarray}
Moreover we remark that there exists $n_3 \equiv
n_3\left(\sigma,s, \kappa\right)$ such that for all $n>n_3$
\begin{eqnarray}
160 \sigma^2 b_n' \log k_n' \leq (200 \sigma)^2 (\gamma)^{1/2} (
\log_2 n)^{7/2}(\log n)^{3/2} .
\end{eqnarray}
We get \ref{prosuperdeltar}, taking $n_0 = n_1 \vee n_2 \vee n_3$.
Similar computations give the result for $\delta_{r',r'}'$.
\end{Pre}

\subsubsection{Proof for the property \ref{superMsup}}

\begin{Lem} \label{profondamsup} There exists $h>0$ such that if \ref{hyp1}, \ref{hyp0} hold and
for all $\kappa \in]0,\kappa^+[$ \ref{hyp4} holds, for all
$\gamma>0$
 there exists $n_0 \equiv n_0(\gamma,\sigma ,\E\left[|\epsilon_0|^3\right])$ such that for all $n>n_0$
\begin{eqnarray}
& & Q\left[\tM_{>} \geq \tmo+L_n \right] \leq h \left(\frac{
\log_3 n}{ \log_2 n}\right)^{1/2}
 , \label{4eqeq} \\
& & Q\left[\tM_{<} \leq \tmo-L_n \right] \leq h \left(\frac{
\log_3 n}{ \log_2 n}\right)^{1/2}  , \label{4eqeq1}
\end{eqnarray}
 see \ref{Msup} for the definitions of $\tM_{<}$ and $\tM_{>}$ and Definition \ref{super} for $L_n$ one.
\end{Lem}

\begin{Pre}
Denote $f(n)= (\log ( q_n (\log n)^{\gamma})) / (\log n)$, where
$q_n$ is given at the end od Definition \ref{super},  we have
\begin{eqnarray}
  Q\left[\tM_{>} \geq \tmo+L_n \right]
 & \equiv & Q\left[
\inf\left\{m>
\tmo,\ S^n_{m}-S^n_{\tmo} \geq f(n) \right\} \geq \tmo+L_n\right] \\
&=& Q\left[ \inf\left\{m> \tmo,\ |S^n_{m}-S^n_{\tmo}| \geq f(n)
\right\} \geq \tmo+L_n \right] \label{4eq392cc} ,
\end{eqnarray}
because $ \tmo $ is a minimizer of the valley
$\{\tM_0',\tmo,\tM_0\}$ and by definition $\tM_0 \geq M_{>}$.
Using Proposition \ref{8eq18}, we know that there exists $n_1
\equiv n_1\left(\sigma, \E\left[ \left|\epsilon_0
\right|^3\right]\right) $ such that for all $n>n_1$
\begin{eqnarray}
Q\left[ -(\sigma^{-1} \log n)^2 \log_2 n \leq \tmo \leq
(\sigma^{-1}\log n)^2 \log_2 n\right] \geq 1-h \left(\frac{ \log_3
n}{ \log_2 n}\right)^{1/2}
  ,
\end{eqnarray}
 so for all $n>n_1$
\begin{eqnarray}
   & & Q\left[ \inf\left\{m> \tmo,\ |S^n_{m}-S^n_{\tmo}| \geq
f(n)
\right\}  \geq \tmo+L_n\right] \\
 & \leq & \sum_{k=-[(\sigma^{-1} \log n)^2 \log_2
n]-1}^{[(\sigma^{-1}\log n)^2 \log_2 n]+1} Q\left[ \inf\left\{m>
k,\ |S^n_{m}-S^n_{k}| \geq f(n) \right\} \geq k+L_n \right] + h
\left(\frac{ \log_3  n}{ \log_2 n}\right)^{1/2} .
\end{eqnarray}
We get that for all $n>n_1$
\begin{eqnarray}
 &  &  Q\left[ \inf\left\{m> \tmo,\ |S^n_{m}-S^n_{\tmo}|
\geq f(n) \right\} \geq
\tmo+L_n \right] \nonumber \\
 & \leq & 2([(\sigma^{-1}\log n)^2 \log_2 n]+1)
Q\left[U^-_{f(n)} \wedge U^+_{f(n)} \geq L_n \right]  + h
\left(\frac{ \log_3  n}{ \log_2 n}\right)^{1/2}
 \label{4eq398cc} .
\end{eqnarray}
Applying inequality \ref{lem101eq1} we get that there exists $n_2
\equiv n_2\left(\sigma, \E\left[ \left|\epsilon_0
\right|^3\right]\right)$ such that for all $n>n_2$

\begin{eqnarray}
& & Q\left[U^-_{f(n)} \wedge U^+_{f(n)} \geq L_n \right] =
\Er\left(\frac{1}{\log n}\right) .
\end{eqnarray}
Replacing this in \ref{4eq398cc} and using \ref{4eq392cc}, we get
\ref{4eqeq} taking $n_0=n_1 \vee n_2$. The proof of \ref{4eqeq1}
is similar.
\end{Pre}

\subsubsection{Proof of Proposition \ref{profonda}}

We only have to collect the results of the Lemmata \ref{moexiste},
\ref{lem100} and \ref{profondamsup}, of the Propositions
\ref{8eq18}, \ref{profondaeta}, \ref{profondadelta},
\ref{profondamu} and \ref{profondaor} and of the Corollary
\ref{superrb}.

\section{ { \large Standard results on sums of i.i.d. random variables
}}

 We recall that for all $\kappa \in ]0,\kappa^+[$, $ C\equiv C(\kappa)= \E_Q\left[e^{\kappa
\epsilon_0}\right] \vee \E_Q\left[e^{-\kappa \epsilon_0}\right]<
+\infty \label{C}$.

\noindent In this section we recall some elementary results on
sums of i.i.d. random variables satisfying the three hypothesis
\ref{hyp1bb}, \ref{hyp0} and \ref{hyp4}.
 \noindent  We will always work on the right of the
origin, that means with $(S_m,m\in \N)$, by symmetry we obtain the
same results for $m \in \Z_-$.

\noindent \\
The following lemma is an immediate consequence of Bernstein
inequality (see \cite{Renyi}).

\begin{Lem} \label{LemBern} Assume \ref{hyp1bb}, \ref{hyp0} hold and for all $\kappa \in ]0,\kappa^+[$ \ref{hyp4} holds.  For all $q>0$ and $p>0$ such that $q<(\sigma^2p) \wedge \left(\sigma^4 p/(2 C)\right)$ we have
\begin{eqnarray}
Q\left[|S_p|>q\right] \leq 2\exp\left\{-\frac{q^2}{2 \sigma^2
p}\left(1-\frac{2qC}{\sigma^4p}\right)\right\} , \label{bernlem1}
\end{eqnarray}
For all $p>1$, $s>0$ and  $k>1$ such that $\log k < (1+s)32
\sigma^2 p  $, for all $0 \leq j \leq p $ we have
\begin{eqnarray}
Q\left[\left|S_p-S_j\right| > \left(32(1+s)  \sigma ^2 p \log k
\right)^{1/2} \right] \leq 2\exp\left\{-\log k +\frac{(p-j)\log
k}{(1+s)64p}+\frac{(p-j)(\log
k)^{3/2}C}{((1+s)32\sigma^2p)^{3/2}}\right\} . \label{bernlem2}
\end{eqnarray}
\end{Lem}

\noindent \\ The following lemma gives an upper bound to the
largest fluctuation of the potential $(S_r,r \in \R)$ in a block
of length $B$ of a given interval.


\begin{Lem} \label{lem102} Assume \ref{hyp1bb}, \ref{hyp0} hold and for all $\kappa \in ]0,\kappa^+[$ \ref{hyp4} holds. For all $s>0$, all integers $K>1 $ and $B>1$ such that $\log K <
\sigma^2 \kappa^2  B $  we have
\begin{eqnarray}
& &  Q\left[\max_{-K-1 \leq i \leq K }\max_{iB \leq j \leq  (i+1)B
}\max_{iB \leq l \leq  (i+1)B }\left(\left|S_l-S_j\right| \right)>
((1+s)32 \sigma^2 B \log K)^{1/2} \right] \nonumber \\
& \leq & 2 K^{-(s-\Er\left( (\log
K)/B\right)^{1/2})}\left(1+\Er\left(H_{K,B}\right) \right) \ .
\label{lem102eq1}
\end{eqnarray}
where $H_{K,B}=K^{-(1-1/64-\Er\left((\log K) / B\right)^{1/2})}$.
For all $L>1$, $K>1$,  all integers $B>1$ such that $[L]+1=K
\times B $ and all $s>0$ such that $  \log K < (1+s)32 \sigma^2
\sigma^2 \kappa^2 B $, we have
\begin{eqnarray}
& &  Q\left[\max_{-[L]-1 \leq m \leq [L]+1 }\max_{m \leq l \leq
m+B }\max_{m \leq j \leq m+B }\left(\left|S_l-S_j\right| \right)>
((1+s)32
\sigma^2 B \log K)^{1/2} \right] \nonumber \\
& \leq & 2(B+1) K^{-(s-\Er\left( (\log
K)/B\right)^{1/2})}\left(1+\Er\left(H_{K,B}\right) \right)  .
\label{lem102eq2}
\end{eqnarray}
\end{Lem}
\begin{Pre}
 Let us prove \ref{lem102eq1}, let $s>0$, $K>1$ and  $B>1$ two positive integers, denoting
 $q=((1+s)32 \sigma^2 B \log K)^{1/2}$. Using the fact that $(\alpha_i,i\in \Z)$
 are i.i.d. we get
\begin{eqnarray}
& & Q\left[\max_{-K-1 \leq i \leq K }\max_{-iB \leq j \leq  (i+1)B
}\max_{iB \leq l \leq  (i+1)B }\left(\left|S_l-S_j\right|
\right)>q \right]   \leq  1-\left(1-Q\left[2\max_{1 \leq j \leq B
}\left(\left|S_j\right| \right)>q \right]\right)^{2K+2} .
\label{Brei2185}
\end{eqnarray}
By Ottaviani inequality (see for example \cite{Breiman} page 45)
\begin{eqnarray}
Q\left[2\max_{1 \leq j \leq B }\left(\left|S_j\right| \right)>q
\right]\leq \frac{Q\left[\left|S_B\right| >q/4 \right]}{1-\sup_{1
\leq j\leq B}\left(Q\left[\left|S_B-S_j\right|
>q/4 \right]\right)} . \label{brei0}
\end{eqnarray}
 Using
\ref{bernlem1},
we have
\begin{eqnarray}
Q\left[\left|S_B\right| >q/4 \right] \leq 2\exp\left\{-\log
K\left(1+s-\Er\left((\log K)/B\right)^{1/2}\right)\right\} .
\label{brei1}
\end{eqnarray}
Similarly, using \ref{bernlem2}, for all $K>1$ such that $\log K <
(1+s)32 \sigma^2 \kappa ^2 B  $, we have
\begin{eqnarray}
\sup_{0 \leq j \leq B}Q\left[\left|S_B-S_j\right| > q \right] &
\leq & 2 K^{-(1-1/64-\Er\left((\log K)/B\right)^{1/2})} .
\label{eqquipu3}
\end{eqnarray}
\noindent Therefore, inserting \ref{brei1} and \ref{eqquipu3} in
\ref{brei0} we get for all $K>1$ such that $ \log K < (1+s) 32
\sigma^2 \kappa^2  B $
\begin{eqnarray}
\qquad Q\left[2\max_{1 \leq j \leq B }\left(\left|S_j\right|
\right)>((1+s)32 \sigma^2 B \log K)^{1/2} \right] &\leq& 2
K^{-(1+s-\Er\left((\log K)/
B\right)^{1/2})}\left(1+\Er\left(H_{K,B}\right) \right) .
\label{5.4}
 \end{eqnarray}
where $H_{K,B}=K^{-(1-1/64-\Er\left(\log K / B\right)^{1/2})}$.
Inserting \ref{5.4} in \ref{Brei2185} and noticing that $(1-x)^a
\geq 1-ax$ for all $ 0 \leq x \leq 1$ and  $a\geq 1$
 we get \ref{lem102eq1}. \\
\noindent Now we prove \ref{lem102eq2}, let $L>1$,
 $B>1$ an integer and  $K>1$ such that $[L]+1=K\times B$, we have $[K] \times B
\leq [L]+1 \leq ([K]+1) \times B $, we remark that
\begin{eqnarray}
& & \max_{-[L]-1 \leq m \leq [L]+1 }\max_{m \leq l \leq m+B
}\max_{m \leq j
\leq m+B }\left(\left|S_l-S_j\right| \right) \\
& \leq & \max_{0 \leq q \leq B } \max_{-[K]-1 \leq i \leq [K]-1
}\max_{iB+q \leq l \leq (i+1)B+q}\max_{iB+q \leq j \leq
(i+1)B+q}\left(\left|S_l-S_j\right| \right) ,
\end{eqnarray}
therefore we have
\begin{eqnarray}
& & Q\left[\max_{-L \leq m \leq L }\max_{m \leq l \leq m+B
}\max_{m \leq j \leq m+B }\left(\left|S_l-S_j\right| \right)>
((1+s)32 \sigma^2 B \log K)^{1/2} \right] \\
& \leq &  (B+1) \times Q\left[  \max_{-[K]-1 \leq i \leq [K]-1
}\max_{iB \leq l \leq (i+1)B } \max_{iB \leq j \leq
(i+1)B }\left(\left|S_l-S_j\right| \right)> \right. \nonumber \\
&& \left. ((1+s)32 \sigma^2 B \log K)^{1/2} \right].
\end{eqnarray}
Using \ref{lem102eq1} we obtain \ref{lem102eq2}.
\end{Pre}
\\

\begin{Lem}  Assume that for all $\kappa \in ]0,k^+[$ \ref{hyp4} holds, for all integer $L>0$ and all $D>0$  we have   \label{lem100}
\begin{eqnarray}
&&  Q\left[\max_{-L \leq i \leq L }\left(\beta_i/\alpha_i\right)
\leq D^{6/\kappa}\right] \geq 1- D^{-6}(2L+1) \E_Q\left[e^{\kappa
\epsilon_0}\right] \ , \label{lem100eq2}
\\
&&  Q\left[\max_{-L \leq i \leq L }\left(\alpha_i/\beta_i\right)
\leq D^{6/\kappa}\right] \geq 1-D^{-6}(2L+1) \E_Q\left[e^{-\kappa
\epsilon_0}\right] \label{lem100eq4} \ ,
\end{eqnarray}
moreover if $D>2^{1+\kappa/6}$
\begin{eqnarray}
& & Q\left[\max_{-L \leq i \leq L }\left(1/\alpha_i\right) \leq
D^{6/\kappa}\right] \geq 1- D^{-6}
2^{\kappa}(2L+1)\E_Q\left[e^{\kappa \epsilon_0}\right]\ ,
\label{lem100eq1}\\
& & Q\left[\max_{-L \leq i \leq L }\left(1/\beta_i\right) \leq
D^{6/\kappa}\right]  \geq 1- D^{-6}
2^{\kappa}(2L+1)\E_Q\left[e^{-\kappa \epsilon_0}\right] \ .
\label{lem100eq3}
\end{eqnarray}
\end{Lem}
\begin{Pre}
This lemma is a simple consequence of the fact that the random
variables $(\alpha_i,\ i \in \Z)$ are i.i.d.
\end{Pre}

\noindent \\ Recalling  \ref{4.3} and \ref{4.3bbc}, we have :

\begin{Lem} \label{lem101}  Assume \ref{hyp1}, \ref{hyp0}, and \ref{hyp4}. Let $\kappa \in
]0,k^+[$, $a>0$, $c>0$ and let us denote $d=a \vee c $.  There
exists $n_0 \equiv n_0\left(\sigma, \E\left[ \left|\epsilon_0
\right|^3\right]\right) $ such that for all $n>n_0$, $L >
\frac{(2(d \log n))^2}{\sigma^2}+1 $ and $D>1$ we have
\begin{eqnarray}
& & Q\left[U_{a}^{-} \wedge U_{c}^{+} > L \right] \leq 2
q_1^{\frac{L \sigma^2}{(2(d \log n))^2+\sigma^2}} ,
\label{lem101eq1} \\
& & Q\left[U_{a}^{-} < U_{c}^{+} \right] \leq
\frac{1}{c+a}\left(c+ \frac{H_d }{\log n} \right) , \label{lem101eq2} \\
& & Q\left[U_{a}^{-} > U_{c}^{+} \right] \leq
\frac{1}{c+a}\left(a+\frac{H_d}{\log n}\right) . \label{lem101eq3}
\end{eqnarray}
where $q_1=0.7+\frac{3,75 \E_Q\left[|\epsilon_0|^{3}\right]}{(d
\log n) \sigma^{2}}<1$ and $H_{d}=(q_1^{\frac{1}{2}\frac{L
\sigma^2}{(2(d \log n))^2+\sigma^2}})/(1-q_1)+(6\log D) /\kappa+
(L^{3/2}
(C)^{1/2} \sigma) /D^{3} $. \\
\end{Lem}
\begin{Pre} We have
\begin{eqnarray}
Q\left[U_{a}^{-} \wedge U_{c}^{+} > L \right] & \leq &
Q\left[U_{d}^{-} \wedge U_{d}^{+}
> L \right]  =  Q\left[\max_{0 \leq l \leq L}\left|S_l\right|<(d \log
n)\right] .
\end{eqnarray}
Let $b=\left[\frac{(2(d \log n))^2}{\sigma^{2}}\right]+1$, for all
$L>b$ there exists $k \equiv k(b,L)$ such that $ k \times b \leq L
\leq b \times (k+1)$, let us denote $[k]$ the integer  part of
$k$, we easily get that

\begin{eqnarray}
 Q\left[U_{a}^{-} \wedge U_{c}^{+} > L \right] &\leq & \left(Q\left[ \left|\frac{S_b}{\sigma
b^{1/2}}\right|< \frac{2(d \log n)}{\sigma b^{1/2}}
\right]\right)^{[k]} . \label{2eq78b}
\end{eqnarray}
\label{ Berry_Essen} Now we use the Berry-Essen theorem (see
\cite{ChoTei} page 299),  we get
\begin{eqnarray}
Q\left[ \left|\frac{S_b}{\sigma b^{1/2}}\right|< \frac{2(d \log
n)}{\sigma b^{1/2}} \right] \leq 2\int^{1}_{0}
\frac{e^{-x^2}}{\sqrt{2\pi}} dx+\frac{3,75
\E_Q\left[|\epsilon_0|^{3}\right]}{(d \log n)\sigma^{2}} .
\label{2eq80}
\end{eqnarray}
Moreover $2\int^{1}_{0} \frac{e^{-x^2}}{\sqrt{2\pi}} dx<0.7 $,
therefore, using \ref{2eq78b} and \ref{2eq80} we get \ref{lem101eq1}. \\
\noindent
 To prove \ref{lem101eq2} we use
 Wald's identity (see \cite{Neveu}) for the martingale $(S^n_t,\ t
\in \R)$ and the regular stopping time $U=U_{a}^{-} \wedge
U_{c}^{+}$. Using  that  $\E_Q \left[S_U^n\right]=0$ and
$\E_Q\left[\left(S_{U_{a}^{-}}^n+a\right) \un_{U_{a}^{-} <
U_{c}^{+}} \right] \leq 0 $ we get that
\begin{eqnarray}
Q\left[U_{a}^{-}<U_{c}^{+}\right] & \leq &
\frac{c}{c+a}+\frac{1}{c+a}\E_Q\left[(S^n_{U_{c}^{+}}-c)\un_{U_{c}^{+}
\leq U_{a}^{-}}\right] \label{2eq88} .
\end{eqnarray}
We have
\begin{eqnarray}
& & \E_Q\left[(S^n_{U_{c}^{+}}-c)\un_{U_{c}^{+} \leq
U_{a}^{-}}\right]=\E_Q\left[(S^n_{U_{c}^{+}}-c)\un_{U_{c}^{+} \leq
U_{a}^{-},U \geq [L]+1
}\right]+\E_Q\left[(S^n_{U_{c}^{+}}-c)\un_{U_{c}^{+} \leq
U_{a}^{-},U < [L]+1 }\right] . \label{eq2.140}
\end{eqnarray}
For the second term on the right hand side of \ref{eq2.140},
noticing that $ (S^n_i-c)\un_{U_{c}^{+} \leq U_{a}^{-},U=i} \leq
\frac{\epsilon_i}{\log n} \un_{U_{c}^{+} \leq U_{a}^{-},U=i} $ we
have
\begin{eqnarray}
 \E_Q\left[(S^n_{U_{c}^{+}}-c)\un_{U_{c}^{+} \leq U_{a}^{-},U
< [L]+1 }\right]  & \leq & \frac{1}{\log n} \sum_{i=1}^{[L
]}\E_Q\left[(\epsilon_i)\un_{U_{c}^{+} \leq U_{a}^{-},U=i }\right]
. \label{eq2.143}
\end{eqnarray}
 For all $D>1$, we have
\begin{eqnarray}
 \frac{1}{\log n}
\sum_{i=1}^{[L]}\E_Q\left[(\epsilon_i)\un_{U_{c}^{+} \leq
U_{a}^{-},U=i }\right]  &=& \frac{1}{\log
n}\sum_{i=1}^{[L]}\E_Q\left[(\epsilon_i)\un_{U_{c}^{+} \leq
U_{a}^{-},U=i, \max_{1 \leq j\leq  [L]
}\left(\epsilon_j\right)\leq
\frac{6}{\kappa}\log D}\right] \label{eq2.148}\\
&+&\frac{1}{\log n}
\sum_{i=1}^{[L]}\E_Q\left[(\epsilon_i)\un_{U_{c}^{+} \leq
U_{a}^{-},U=i, \max_{1 \leq j\leq  [L] }\left(\epsilon_j\right)>
\frac{6}{\kappa}\log D}\right] \label{eq2.149} \\
& \leq & \frac{6\log D}{\kappa \log n}+\frac{ \sigma [L]}{\log n}
\left(Q\left[\max_{1 \leq j\leq  [L] }\left(\epsilon_j\right)>
\frac{6}{\kappa}\log D\right] \right)^{1/2} , \label{4eq89}
\end{eqnarray}
where we have used that for the sum in the right hand side of
\ref{eq2.148} the $\epsilon_i$ are bounded by
$\frac{6}{\kappa}\log D$ and for the sum \ref{eq2.149} the
Cauchy-Schwarz inequality.
 To end we use \ref{lem100eq2}, for all $D>2^{1+\kappa/6}$
\begin{eqnarray}
\E_Q\left[(S^n_{U_{c}^{+}}-c)\un_{U_{c}^{+} \leq U_{a}^{-},U <
[L]+1 }\right]&\leq & \frac{6\log D}{ \kappa \log n}+\frac{ \sigma
([L])^{3/2} \left(\E_Q\left[e^{\kappa \log
\left(\frac{\beta_0}{\alpha_0}\right)}\right]\right)^{1/2} }{D^{3}
\log n }. \label{2eq97}
\end{eqnarray}
For the first term of the right hand side of \ref{eq2.140}, using
 Cauchy-Schwarz inequality we get
\begin{eqnarray}
\E_Q\left[(S^n_{U_{c}^{+}}-c)\un_{U_{c}^{+} \leq U_{a}^{-},U \geq
[L]+1 }\right] & \leq & \frac{\sigma}{\log n}
\sum_{i=[L]+1}^{\infty} \left(Q\left[U \geq i
\right]\right)^{1/2},
\end{eqnarray}
then, to estimate, $Q\left[U \geq i \right]$ we use
\ref{lem101eq1}.
Collecting what we did above we get \ref{lem101eq2}.
\end{Pre}

 \noindent \\ We use the following notation  $Q[.|S_0=y]=Q_y[.]$ ($Q[.|S_0=0] \equiv Q_0[.] = Q[.]$).

\begin{Rem}
 $\bullet$ For all $a>0$, $b>0$ and $l>0$ we have
\begin{eqnarray}
Q\left[U^+_c > l\right] \leq Q\left[U^+_c \wedge U^-_a > l\right]+
Q\left[U^+_c > U^-_a \right] \label{symb} \label{sym} .
\end{eqnarray}
\end{Rem}

\newpage

\begin{figure}[h]
\begin{center}
\input{thfig4.pstex_t} \caption{} \label{thfignew1}
\input{thfig3.pstex_t} \caption{} \label{thfig3}
\end{center}
\end{figure}

\noindent \\ \textbf{ Acknowledgment  } This article is a part of
my Phd thesis made under the supervision of P. Picco. I would like
to thank him for helpful discussions all along the last three
years. I would like to thank R. Correa and all the members of the
C.M.M. (Santiago, Chili) for their hospitality during all the year
2002.

{\small \bibliography{article}}

\begin{thebibliography}{49}
\expandafter\ifx\csname natexlab\endcsname\relax\def\natexlab#1{#1}\fi
\expandafter\ifx\csname url\endcsname\relax
  \def\url#1{{\tt #1}}\fi

\bibitem[Solomon(1975)]{Solomon}
F.~Solomon.
\newblock Random walks in random environment.
\newblock {\em Ann. Probab.}, \textbf{3}\penalty0 (1):\penalty0 \ 1--31, 1975.

\bibitem[Kesten et~al.(1975)Kesten, Kozlov, and Spitzer]{KesKozSpi}
H.~Kesten, M.V. Kozlov, and F.~Spitzer.
\newblock A limit law for random walk in a random environment.
\newblock {\em Comp. Math.}, \textbf{30}:\penalty0 \ 145--168, 1975.

\bibitem[Sinai(1982)]{Sinai}
Ya.~G. Sinai.
\newblock The limit behaviour of a one-dimensional random walk in a random
  medium.
\newblock {\em Theory Probab. Appl.}, \textbf{27}\penalty0 (2):\penalty0 \
  256--268, 1982.

\bibitem[Golosov(1984)]{Golosov}
A.~O. Golosov.
\newblock Localization of random walks in one-dimensional random environments.
\newblock {\em Communications in Mathematical Physics}, \textbf{92}:\penalty0 \
  491--506, 1984.

\bibitem[Golosov(1986)]{Golosov0}
A.~O. Golosov.
\newblock Limit distributions for random walks in random environments.
\newblock {\em Soviet Math. Dokl.}, \textbf{28}:\penalty0 \ 18--22, 1986.

\bibitem[Kesten(1986)]{Kesten2}
H.~Kesten.
\newblock The limit distribution of {S}inai's random walk in random
  environment.
\newblock {\em Physica}, \textbf{138A}:\penalty0 \ 299--309, 1986.

\bibitem[Deheuvels and Révész(1986)]{Deh&Revesz}
P.~Deheuvels and P.~Révész.
\newblock Simple random walk on the line in random environment.
\newblock {\em Probab. Theory Related Fields}, \textbf{72}:\penalty0 \
  215--230, 1986.

\bibitem[R\'ev\'esz(1989)]{Revesz}
P.~R\'ev\'esz.
\newblock {\em Random walk in random and non-random environments}.
\newblock World Scientific, 1989.

\bibitem[Greven and Hollander(1994)]{GreHol2}
A.~Greven and F.~Hollander.
\newblock Large deviation for a walk in random environment.
\newblock {\em Ann. probab.}, \textbf{27}\penalty0 (4):\penalty0 \ 1381--1428,
  1994.

\bibitem[Zeitouni and Gantert(1998)]{ZeiGan}
O.~Zeitouni and N.~Gantert.
\newblock Quenched sub-exponential tail estimates for one-dimentional random
  walk in random environment.
\newblock {\em Comm. Math. Phys.}, \textbf{194}:\penalty0 \ 177--190, 1998.

\bibitem[Pisztora and Povel(1999)]{PiPO}
A.~Pisztora and T.~Povel.
\newblock Large deviation principle for random walk in a quenched random
  environment in the low speed regime.
\newblock {\em Ann. Probab.}, \textbf{27}:\penalty0 \ 1389--1413, 1999.

\bibitem[Zeitouni et~al.(1999)Zeitouni, Pisztora, and Povel]{PiPoZe}
O.~Zeitouni, A.~Pisztora, and T.~Povel.
\newblock Precise large deviation estimates for a one-dimensional random walk
  in a random environment.
\newblock {\em Probab. Theory Related Fields}, \textbf{113}:\penalty0 \
  191--219, 1999.

\bibitem[Comets et~al.(2000)Comets, Zeitouni, and Gantert]{ZeCoGa}
F.~Comets, O.~Zeitouni, and N.~Gantert.
\newblock Quenched, annealed and functional large deviations for
  one-dimensional random walk in random environment.
\newblock {\em Probab. Theory Related Fields}, 118,:\penalty0 65--114, 2000.

\bibitem[Zeitouni(2001)]{Zeitouni}
O.~Zeitouni.
\newblock Lectures notes on random walks in random environment.
\newblock {\em St Flour Summer School}, 2001.

\bibitem[Shi(1998)]{Shi}
Z.~Shi.
\newblock A local time curiosity in random environment.
\newblock {\em Stochastic Process. Appl.}, \textbf{76}\penalty0 (2):\penalty0 \
  231--250, 1998.

\bibitem[Hu and Shi(1998{\natexlab{a}})]{HuShi2}
Y.~Hu and Z.~Shi.
\newblock The limits of {S}inai's simple random walk in random environment.
\newblock {\em Ann. Probab.}, \textbf{26}\penalty0 (4):\penalty0 \ 1477--1521,
  1998{\natexlab{a}}.

\bibitem[Hu and Shi(1998{\natexlab{b}})]{HuShi1}
Y.~Hu and Z.~Shi.
\newblock The local time of simple random walk in random environment.
\newblock {\em J. of Theoret. Probab.}, \textbf{11}\penalty0 (3),
  1998{\natexlab{b}}.

\bibitem[Hu(2000{\natexlab{a}})]{Hu1}
Y.~Hu.
\newblock The logarithmic average of {S}inai's walk in random environment.
\newblock {\em Period. Math. Hungar.}, \textbf{41}:\penalty0 175--185,
  2000{\natexlab{a}}.

\bibitem[Hu(2000{\natexlab{b}})]{Hu}
Y.~Hu.
\newblock Tightness of localization and return time in random environment.
\newblock {\em Stochastic Process. Appl.}, \textbf{86}\penalty0 (1):\penalty0 \
  81--101, 2000{\natexlab{b}}.

\bibitem[Hu and Shi(2000)]{HuShi0}
Y.~Hu and Z.~Shi.
\newblock The problem of the most visited site in random environment.
\newblock {\em Probab. Theory Related Fields}, \textbf{116}\penalty0
  (2):\penalty0 \ 273--302, 2000.

\bibitem[Schumacher(1985)]{Schumacher}
S.~Schumacher.
\newblock Diffusions with random coefficients.
\newblock {\em Contemp. Math.}, \textbf{41}:\penalty0 \ 351--356, 1985.

\bibitem[Brox(1986)]{Brox}
T.~Brox.
\newblock A one-dimensional diffusion process in a {W}iener medium.
\newblock {\em Ann. Probab.}, \textbf{14}\penalty0 (4):\penalty0 \ 1206--1218,
  1986.

\bibitem[Shi(2001)]{Shi1}
Z.~Shi.
\newblock {S}inai's walk via stochastic calculus.
\newblock {\em Panoramas et Synthèses}, \textbf{12}:\penalty0 \ 53--74, 2001.

\bibitem[Dembo et~al.(2001)Dembo, Guionnet, and Zeitouni]{ZeDEGu}
A.~Dembo, A.~Guionnet, and O.~Zeitouni.
\newblock Aging properties of {S}inai's model of random walk in random
  environment.
\newblock {\em In St. Flour summer school 2001 lecture notes}, 2001.

\bibitem[Comets and Popov(2003)]{ComPop}
F.~Comets and S.~Popov.
\newblock Limit law for transition probabilities and moderate deviations for
  {S}inai's random walk in random environment.
\newblock {\em Preprint}, 2003.

\bibitem[Gantert and Shi(2002)]{GanShi}
N.~Gantert and Z.~Shi.
\newblock Many visits to a single site by a transient random walk in random
  environment.
\newblock {\em Stochastic Processes and their applications}, 99:\penalty0 \
  159--176, 2002.

\bibitem[Tanaka(1994)]{Tanaka}
H.~Tanaka.
\newblock Localization of a diffusion process in a one-dimensional brownian
  environmement.
\newblock {\em Comm. Pure Appl. Math.}, \textbf{17}:\penalty0 \ 755--766, 1994.

\bibitem[Mathieu(1995)]{Mathieu2}
P.~Mathieu.
\newblock Limit theorems for diffusions with a random potential.
\newblock {\em Stochastic Process. Appl.}, \textbf{60}:\penalty0 \ 103--111,
  1995.

\bibitem[Tanaka(1997)]{Tanaka2}
H.~Tanaka.
\newblock Limit theorem for a brownian motion with drift in a white noise
  environment.
\newblock {\em Chaos Solitons Fractals}, \textbf{11}:\penalty0 \ 1807--1816,
  1997.

\bibitem[Tanaka and Kawazu(1997)]{KawTan}
H.~Tanaka and K.~Kawazu.
\newblock A diffusion process in a brownian environment with drift.
\newblock {\em J. Math. Soc. Japan}, \textbf{49}:\penalty0 \ 189--211, 1997.

\bibitem[Mathieu(1998)]{Mathieu1}
P.~Mathieu.
\newblock On random perturbations of dynamical systems and diffusion with a
  random potentiel in dimension one.
\newblock {\em Stochastic Process. Appl.}, \textbf{77}:\penalty0 \ 53--67,
  1998.

\bibitem[Taleb(2001)]{Taleb0}
M.~Taleb.
\newblock Large deviations for a brownian motion in a drifted brownian
  potential.
\newblock {\em Ann. Probab.}, \textbf{29}\penalty0 (3):\penalty0 \ 1173--1204,
  2001.

\bibitem[Kalikow(1981)]{Kalikow}
S.~A. Kalikow.
\newblock Generalised random walk in random environment.
\newblock {\em The Ann. of Prob.}, \textbf{9}\penalty0 (5):\penalty0 \
  753--768, 1981.

\bibitem[Anshelevich et~al.(1982)Anshelevich, Khanin, and Sinai]{AnKhSi}
V.~V. Anshelevich, K.~M. Khanin, and Ya.~G. Sinai.
\newblock Symmetric random walks in random environments.
\newblock {\em Com. Math. Phy.}, \textbf{85}:\penalty0 \ 449--470, 1982.

\bibitem[Durrett(1986)]{Durett}
R.~Durrett.
\newblock Some multidimensional rwre with subclassical limiting behavior.
\newblock {\em Commun. Math. Phys.}, \textbf{104}:\penalty0 \ 87--102, 1986.

\bibitem[Bouchaud et~al.(1987)Bouchaud, Comtet, Georges, and Doussal]{BoCoGeDo}
J.P. Bouchaud, A.~Comtet, A.~Georges, and P.~Le Doussal.
\newblock Anomalous diffusion in random media of any dimensionality.
\newblock {\em J. Physique}, \textbf{48}:\penalty0 \ 1445--1450, 1987.

\bibitem[Bricmont and Kupiainen(1991)]{BricKup}
J.~Bricmont and A.~Kupiainen.
\newblock Random walks in asymetric random environments.
\newblock {\em Comm. in Math. Phys.}, \textbf{142}:\penalty0 342--420, 1991.

\bibitem[Sznitman(1999)]{Sznitman}
A.~S. Sznitman.
\newblock Lectures on random motions in random media.
\newblock {\em Preprint}, 1999.

\bibitem[Sznitman(2003)]{Sznitman3}
A.~S. Sznitman.
\newblock On new examples of ballistic random walks in random environment.
\newblock {\em Ann. Probab.}, \textbf{31}\penalty0 (1):\penalty0 \ 285--322,
  2003.

\bibitem[Varadhan(2003)]{Varadhan}
S.~R.~S. Varadhan.
\newblock Large deviations for random walks in random environment.
\newblock {\em Comm. Pure Appl. Math.}, \textbf{56}\penalty0 (8):\penalty0 \
  1222--1245, 2003.

\bibitem[Rassoul-Agha(2003)]{Rass2}
F.~Rassoul-Agha.
\newblock The point of view of the particule on the law of large numbers for
  random walks in a mixing random environment.
\newblock {\em Ann. Probab.}, \textbf{31}:\penalty0 \ 1441--1463, 2003.

\bibitem[Comets and Zeitouni(2004)]{ZeiCom2}
F.~Comets and O.~Zeitouni.
\newblock A law of large numbers for random walk in random environments.
\newblock {\em To appear in Ann. Probab.}, 2004.

\bibitem[Chung(1967)]{Chung}
K.~L. Chung.
\newblock {\em Markov Chains}.
\newblock Springer-Verlag, 1967.

\bibitem[Cassandro et~al.(2004+)Cassandro, Orlandi, Picco, and Varés]{Picco}
M.~Cassandro, E.~Orlandi, P.~Picco, and M.~E. Varés.
\newblock One dimensional random field kac's model : localisation of the
  phases.
\newblock {\em Preprint}, 2004+.

\bibitem[LeCam(1986)]{Lecam}
L.~LeCam.
\newblock {\em Asymptotic methods in statistical decision theory}.
\newblock Springer-Verlag, 1986.

\bibitem[Renyi(1970)]{Renyi}
A.~Renyi.
\newblock {\em Probability Theory}.
\newblock North-Holland Publishing Company, 1970.

\bibitem[Breiman(1968)]{Breiman}
L.~Breiman.
\newblock {\em Probability}.
\newblock Addison-Wesley Publishing Company, Inc, 1968.

\bibitem[Chow and Teicher(1997)]{ChoTei}
Y.~S. Chow and H.~Teicher.
\newblock {\em Probability Theory}.
\newblock Srpinger, third edition, 1997.

\bibitem[Neveu(1972)]{Neveu}
J.~Neveu.
\newblock {\em Martinguales à temps discret}.
\newblock Masson et Cie, 1972.

\end{thebibliography}

\vspace{1cm} \noindent
\begin{tabular}{ll}
Centre de Physique Théorique - C.N.R.S.  & \hspace{2.3cm} Centro de Modelamiento Matem\'atico - C.N.R.S. \\
Université Aix-Marseille II & \hspace{2.3cm} Universidad de Chile  \\
Luminy Case 907 & \hspace{2.3cm} Blanco Encalada 2120 piso 7 \\
13288 Marseille cedex 09, France & \hspace{2.3cm} Santiago de
Chile
\end{tabular}

\end{document}